\documentclass[11pt]{amsart}
\usepackage{t1enc}
\usepackage[english]{babel}
\usepackage{amsmath,amsthm,amssymb}
\usepackage{amsfonts} 
\usepackage{latexsym} 
\usepackage[dvips]{graphicx}
\usepackage{graphicx}
\usepackage[natural]{xcolor}
\usepackage{float}
\usepackage{enumerate}
\usepackage[subnum]{cases}
\usepackage{multirow}
\usepackage{hyperref}
\usepackage{psfrag}
\usepackage{bm}
\usepackage{mathrsfs}
\usepackage{graphicx}
\usepackage{comment}
\newcommand{\pf}{\noindent\textbf{Proof}.\quad}
\usepackage{listings}
\newtheorem{definition}{Definition}[section]
\newtheorem{thm}[definition]{Theorem}

\newtheorem{lemma}{Lemma}[section]

\usepackage{lineno}
\newcommand{\D}{\Delta}

\newcommand{\phiv}{\varphi}

\newcommand{\phibar}{\overline{\phiv}}
\newcommand{\pibar}{\overline{\pi}}

\newcommand{\xibar}{\overline{\xi}}

\newcommand{\epv}{\varepsilon}

\newcommand{\dd}{\delta}

  \title [Goldberg-Seymour conjecture]
{A short proof of the Goldberg-Seymour conjecture}
\author{Guantao Chen}
\thanks{Partially supported by NSF grants DMS-1855716 and DMS-2154331.}
\email{gchen@gsu.edu}
\address{Department of Mathematics and Statistics, Georgia State University, Atlanta, GA 30303}
\author{Yanli Hao}
\email{yhao98@gatech.edu}
\address{School of Mathematics, Georgia Institute of Technology, Atlanta, GA 30332}
\author{Xingxing Yu}
\thanks{partially supported by NSF grant DMS-2348702}
\email{yu@math.gatech.edu}
\address{School of Mathematics, Georgia Institute of Technology, Atlanta, GA 30332}
\author{Wenan Zang}
\email{wzang@maths.hku.hk}
\thanks{Supported in part by the Research Grants Council of Hong Kong}
\address{Department of Mathematics, The University of Hong Kong, Hong Kong, China }

\date{}
\begin{document}

\begin{abstract}

For a multigraph $G$, $\chi'(G)$ denotes the chromatic index of $G$, $\Delta(G)$ the maximum degree of $G$, and $\Gamma(G) = \max\{\lceil 2|E(H)|/(|V (H)|-1)\rceil: H \subseteq G \mbox{ and } |V(H)| \mbox{ odd}\}$. 
As a generalization of Vizing’s classical coloring result for simple graphs, the
Goldberg-Seymour conjecture, posed in the 1970s, states that $\chi'(G)=\max\{\Delta(G), \Gamma(G)\}$ or $\chi'(G)=\max\{\Delta(G) + 1, \Gamma(G)\}$. Hochbaum, Nishizeki and Shmoys further conjectured in 1986 that such a coloring can be found in polynomial time. A long proof of the Goldberg-Seymour conjecture was announced in 2019 by Chen, Jing, and Zang, and one case in that proof was eliminated recently by Jing (but the proof is still long); and neither proof has been verified. In this paper, we give a proof of the Goldberg-Seymour conjecture that is significantly shorter and confirm the Hochbaum-Nishizeki-Schmoys conjecture by providing an $O(|V|^5|E|^3)$ time algorithm for finding a $\max\{\Delta(G) + 1, \Gamma(G)\}$-edge-coloring of $G$.

\end{abstract}

\maketitle

\newpage

\setcounter{page}{1}

\section{Introduction}
Let $G$ be a loopless multigraph. Denote by  $V(G)$ and $E(G)$  the vertex set and edge
set of $G$, respectively, $\Delta(G)$ the maximum
degree of $G$, and $\chi'(G)$ the chromatic index of $G$. Vizing’s classical theorem states
that $\chi'(G) \in\{ \Delta(G), \Delta(G)+ 1\}$ when $G$ is simple, which no longer
holds for multigraphs. Shannon~\cite{Shannon49} showed that $\chi'(G) \le \frac 32 \D(G)$; and Vizing~\cite{Vizing68} and Gupta~\cite{Gupta67} showed, independently, that $\chi'(G) \le \D(G) + \mu(G)$, where $\mu(G)$ is the maximum number of parallel edges between vertices of $G$.   Let
\[
\Gamma (G) = \max\{\lceil 2|E(H)|/(|V (H)|-1) \rceil: H\subseteq G \mbox{ and $|V(H)|$ odd} \}.
\]
It is a simple exercise to verify that $\chi'(G) \ge  \Gamma(G)$; so $\chi'(G) \ge  \max\{\Delta(G), \Gamma(G)\}$. Using Edmonds' matching polytope theorem, Seymour~\cite{Seymour74} showed that $\max\{\D(G), \Gamma(G)\}$ can be computed in polynomial time. Chen, Zang and Zhao~\cite{ChenZZ19} recently gave a combinatorial algorithm that computes $\Gamma(G)$ in polynomial time. 

Goldberg \cite{Goldberg73} and Seymour \cite{Seymour74} independently conjectured in 1973 that, for any multigraph
$G$, $\chi'(G)\le \max\{\Delta(G)+1, \Gamma(G)\}$, which would imply that the chromatic index of any
multigraph can be approximated in polynomial time within additive constant 1. Hochbaum, Nishizeki and Shmoys \cite{Hochbaum86} further conjectured in 1986 that a $\max\{\Delta(G)+1, \Gamma(G)\}$-edge-coloring of $G$ can be found in polynomial time. We refer the reader to~\cite{StiebitzBook,CJZ22} for literature on the history and development of these conjectures. We mention in particular that the Goldberg-Seymour conjecture holds asymptotically \cite{Kahn96, CYZ12, Scheide10, ChenGKPS18} and for random graphs \cite{HaxellKK19}. 

A proof of the Goldberg-Seymour conjecture was announced in 2019 \cite{CJZ22} and a more recent proof \cite{jing2024} was announced where a case in \cite{CJZ22} was removed.
Both \cite{CJZ22} and \cite{jing2024} are very long and have not been independently verified. Here we give a significantly
shorter proof of the Goldberg-Seymour conjecture and provide an $O(|V|^5|E|^3)$ time algorithm that produces the desired coloring thereby also establishing the Hochbaum-Nishizeki-Shmoys conjecture.

\begin{thm}\label{thm:main}
For any multigraph $G$, one can find a $\max\{\Delta(G)+1, \Gamma(G)\}$-edge-coloring of $G$ in $O(|V(G)|^5|E(G)|^{3})$ time.
\end{thm}

We made no attempt to optimize the complexity bound, but we think it can be improved. The main idea of our proof of Theorem~\ref{thm:main} is to extend a partial coloring of $G$  through color exchanges based on a special type of trees.  
We will give a brief overview of the proof in Section 2 after introducing the necessary concepts. We will then prove lemmas in Section 3 on color exchanges that preserve the tree structures, and prove lemmas in Section 4 on reduction to smaller trees. We complete the proof of Theorem~\ref{thm:main} in Section 5. 

In the remainder of this section, we introduce the notation for this paper. Let $G$ be a graph. For $X \subseteq V(G)$, we use $G[X]$ to denote the subgraph of $G$ induced by $X$ and let $\partial_G X$ denote the set of edges in $E(G)$ with
exactly one end in $X$.  For the set $S\subseteq E(G)$, we use $G-S$ to denote the graph obtained from $G$ by removing all edges in $S$. Let $H\subseteq G$ (i.e., $H$ is a subgraph of $G$).  We write $G[H]$ for $G[V(H)]$, $\partial_G H$ for $\partial_G V(H)$, and $G-H := G[V(G) \backslash V(H)]$. We drop the subcript $G$ if no confusion arises.
An {\bf $H$-path} is a path in $G$ that has both ends
in $V (H)$ but is otherwise disjoint from $H$. For a set $S\subseteq E(G)$, we use $G[S]$ to denote the subgraph of $G$ induced by $S$ and $H+S$ to denote the graph with vertex set $V(H)\cup V(G[S])$ and the edge set $E(H)\cup S$.  When $S=\{g\}$ we write $H+g$ instead of $H+\{g\}$. 
For any two integers $h, k$ with $h \le k$, let
$[h, k]=\{h,  \ldots, k\}$; let $[k]:=[1,k]$ if $k>0$; and let $[k]=\emptyset$ if $k\le 0$.

We call a tree $T$ an {\it ordered} tree  if its edges are linearly ordered as $e_1<e_2<\cdots < e_m$, such that, for every $i\in [m]$, $T[\{e_1, \ldots, e_i\}]$ is a tree and said to be a {\it prefix} of $T$.  For any edge $e_i\in E(T)$,  let $T(e_i):= T[\{e_1, \dots, e_{i}\}]$; and for any vertex $v\in V(T)$, we use $T(v)$ to denote the minimal tree $T(e_i)$ such that $v\in V(T(e_i))$. Let $T^-:=T[\{e_1, \ldots, e_{m-1}\}]$.

  We will use natural numbers to denote colors. For a graph $G$ and a set $S\subseteq E(G)$, a {\it partial coloring} of $G$  with {\it support} $S$ (or simply {\it coloring}) is a function $\phiv: S \to [k]$ (for some positive integer $k$) such that no two incident edges of $G$ receive the same color.  For any $v\in V(G)$, we write $\phiv(v):=\{\phiv(e): e\in S \mbox{ and $e$ is incident with $v$}\}$ and $\phibar(v)=[k]\setminus \phiv(v)$. For any $U\subseteq V(G)$, let $\phibar(U)=\bigcup_{v\in U}\phibar(v)$. For any subgraph $H$ of $G$, let $\phibar(H)=\phibar(V(H))$ and $\phiv(H)=\phiv(E(H))$. 
  Let $\phibar(g)=\phibar(G[g])$ for $g\in E(G)$. 
  
  Let $G$ be a graph and $\phiv$ a partial coloring of $G$. 
For an ordered tree $T$ in $G$, we grow $T$ by adding edges $e_1, \ldots, e_m$ in order, such that $T+\{e_1, \ldots, e_i\}$ is a tree for each $i\in [m]$,  and $\phiv(e_1)\in \phibar(T)$ and, for $2\le i\le m$,  $\phiv(e_i)\in \phibar(T+\{e_1, \ldots, e_{i-1}\})$. We call a maximal such tree a {\bf closure of $T$ under $\phiv$} and denote it by $cl_\phiv(T)$. Observe that\\ 
(1) $cl_\phiv(T)$ can be found in $O(|V(G)|^2k^2)$ time, and\\
(2) $cl_{\phiv}(T)$ is not unique, but $V(cl_{\phiv}(T))$ is uniquely determined. To see this, let $f_1, \ldots, f_n$ denote another sequence that gives a closure of $T$ and let $i\in [n]$ be minimum such that an incident vertex of $f_i$ is not in $cl_{\phiv}(T)$. Then $\phiv(f_i)\in \phibar(T+\{f_1, \ldots, f_{i-1}\}])\subseteq \phibar(cl_{\phiv}(T))$, contradicting the definition of closure.  
  
   Let $G$ be a graph and $\phiv$ a partial coloring of $G$, with $\phiv(G)\subseteq [k]$. For two colors $\alpha, \beta\in [k]$, each component of $G[\phiv^{-1}(\alpha, \beta)]$ 
  is a path or cycle, and a component that is a path is called an {\bf $\{\alpha, \beta\}$-chain}. For $H\subseteq G$, an $H$-path contained in $G[\phiv^{-1}(\alpha, \beta)]$ is called an {\bf $(H, \phiv,\alpha,\beta)$-path}. Let $P$ be a component of $G[\phiv^{-1}(\alpha, \beta)]$; we use $\phiv/P$ to denote the coloring obtained from $\phiv$ by swapping colors $\alpha$ and $\beta$ on $P$ and we say that $\phiv/P$ is obained from $\phiv$ by a {\it Kempe exchange}.  Two partial colorings of $G$ are {\it Kempe equivalent} if one can be obtained from the other by a sequence of Kempe exchanges. If $H\subseteq G$ and 
 $\alpha,\beta\notin \partial H$, then $\phiv/(G-H,\alpha,\beta)$ denotes the coloring obtained from $\phiv$ by exchanging colors $\alpha, \beta$ on all edges of $G-H$, which is clearly is Kempe equivalent to $\phiv$.

\section{Tashkinov trees}
Let $G$ be a multigraph, $k = \max\{\D(G) +1, \Gamma(G)\}$, $\phiv$ a partial coloring of $G$ using colors from $[k]$, and $g\in E(G)\setminus \phiv^{-1}([k])$. We wish to obtain a coloring $\pi$ Kempe equivalent to $\phiv$, such that $\phiv^{-1}([k])\cup \{g\}\subseteq \pi^{-1}([k])$. We do this by growing a tree in various stages.

We say that a tree $T\subseteq G$ is {\bf $\phiv$-elementary} if
$\phibar(u)\cap \phibar(v)=\emptyset$ for any distinct $u,v\in V(T)$; $T$ is {\bf $\phiv$-closed} if  $\phiv(\partial T)\cap  \phibar(T) = \emptyset$. 
Note that 
if $T$ is both $\phiv$-elementary and $\phiv$-closed then there exists $\alpha\in [k]$ such that $|\partial T\cap \phiv^{-1}(\alpha) |\ge 2$; for, otherwise, for each $\beta\in [k]$, the color class $E(G[T])\cap \phiv^{-1}(\beta)$ is a matching of $G[T]$ of size $\lfloor |V(T)|/2\rfloor$, giving $|E(G[T])| \ge k\cdot (|V(T)|-1)/2 +1$, and hence $\Gamma(G) > k$, a contradiction. In particular, if $V(T)=V(G)$ then $T$ is $\phiv$-closed and, hence, not $\phiv$-elementary. We will construct trees $T$ that are not $\phiv$-elementary, and reduce these trees via Kemple exchanges to get a coloring $\pi$ of $G$ such that $G[g]$ is not $\pi$-elementary, hence extending the coloring $\pi$ by also coloring $g$.

Let $T:=cl_\phiv(G[g])$ and call it a {\bf $(0,0,0)$-tree}. 
If $T$ is not $\phiv$-elementary, then Tashkinov in \cite{Tashkinov20} reduces $T$ to a smaller non-elementary tree containing $g$ under certain lexicographic ordering, and proves the following result. See for example,  \cite{StiebitzBook} (Lemma 5.2) and \cite{Sebo22+}.

\begin{lemma}\label{lem:recoloring}
Let $G$ be a multigtaph, $g\in E(G)$, and $\phiv$ a partial coloring of $G-g$ with $\phiv(G) \subseteq [k]$, where $k = \max\{\D(G) +1, \Gamma(G)\}$. If $cl_{\phiv}(G[g])$ is not $\phiv$-elementary, then there is a coloring $\pi$ of $G$, such that  $\pi^{-1}([k])=\phiv^{-1} ([k])\cup \{g\}$ and $\pi$ is Kempe equivalent to $\phiv$ on $G-g$. Moreover, such $\pi$ can be found in $O(|V(G)|^3)$ time.
\end{lemma}

The main challenge to proving the Goldberg-Seymour conjecture (and the Hochbaum-Nishizeki-Shmoys conjecture) is to deal with the case when $cl_{\phiv}(G[g])$ is $\phiv$-elementary and $\phiv$-closed. An approach, first explored by Chen, Yu and Zang \cite{CYZ12}, is to further grow $cl_{\phiv}(G[g])$ by including certain edges $e\in \partial cl_{\phiv}(G[g])$ with $\phiv(e)\notin \phibar(cl_\phiv(G[g]))$.

Let $T_0=cl_{\phiv}(G[g])$ and suppose it is
$\phiv$-elementary and $\phiv$-closed.
We wish to further
grow $T_0$. Fix a color $\alpha_0\in \phibar(g)$ and choose a color $\delta_0$ such that
$|\partial T_0\cap \phiv^{-1}(\delta_0)|\ge 2$. Let $P$ denote a $(T_0,\phiv,\alpha_0,\dd_0)$-path in $G$ (existence guaranteed by Lemma \ref{thm-2-Tn-extension}) and any path obtained from $P$ by removing one end vertex of $P$ is called a {\bf $(T_0,\phiv,\alpha_0,\delta_0)$-nonexit path}.
Let $T_0^+(\phiv,\alpha_0,\dd_0)$ be an ordered tree obtained from $T_0$ by adding $(T_0,\phiv,\alpha_0,\delta_0)$-nonexit paths, one for each $(T_0,\phiv,\alpha_0,\dd_0)$-path,  so that those paths are added in an arbitrary order but the edges in each such path are added consecutively in the order along that path starting from the edge incident with $T_0$. We then grow $T_0^+(\phiv,\alpha_0,\dd_0)$ to $T_1:=cl_{\phiv}(T_0^+(\phiv,\alpha_0,\dd_0))$.
If $T_1$ is not $\phiv$-elementary, we will obtain a `smaller' tree that is not $\pi$-elementary for some coloring $\pi$ Kempe equivalent to $\phiv$ (see Section 4).
If $T_1$ is $\phiv$-elementary and $\phiv$-closed, 
we choose a color $\delta_1\in [k]$ with $|\partial T_1\cap \phiv^{-1}(\dd_1)| \ge 2$, and say that $T_1$ is a (1,0,0)-tree with {\bf base} $T_0$, {\bf stage} $T_1$,  and ordered multiset of {\bf connecting colors} $\{\delta_0,\dd_1\}$. We further grow $T_1$ to $T_1^+(\phiv,\alpha_0,\dd_1)$ by adding  $(T_1,\phiv,\alpha_0,\delta_1)$-nonexit paths (guaranteed to exist by Lemma \ref{thm-2-Tn-extension}), and set $T_2=cl_{\phiv}(T_1^+(\phiv,\alpha_0,\dd_1))$ which has base $T_0$, stages $T_1,T_2$. We give a formal definition below.

We use {\bf quadruple} $(G, g,\phiv,\alpha_0)$ to denote a graph $G$ with $g\in E(G)$, a partial coloring $\phiv$ of $G-g$, $\alpha_0\in \phibar(g)$,  and $\phiv(G)\subseteq [k]$, where $k=\max\{\Delta(G)+1,\Gamma(G)\}$. Throughout the remainder of this section, we fix a quadruple $(G, g,\phiv,\alpha_0)$.

\begin{definition}\label{defi-3-TP}
Let $n$ be a positive integer. Suppose there exist ordered trees $T_0, \ldots, T_{n}$ with $T_0=cl_{\phiv}(G[g])$ and for $i\in [0,n]$ there exist $\delta_i\in [k]$ with $|\partial T_i\cap \phiv^{-1}(\dd_i)|\ge 2$, such that,  for $i\in [n]$, $T_{i}=cl_{\phiv}(T_{i-1}^+(\phiv,\alpha_0,\dd_{i-1}))$, where $T_{i-1}^+(\phiv,\alpha_0,\dd_i)$ is an ordered tree obtained from $T_{i-1}$ by adding $(T_{i-1},\phiv,\alpha_0,\delta_i)$-nonexit paths, {\bf one} for each $(T_{i-1},\phiv,\alpha_0,\dd_i)$-path by removing one end vertex, so that those paths are added in an arbitrary order but the edges in each such path are added consecutively in the order along that path starting from the edge incident with $T_0$. We say that $T_{n}$ is an {\bf (n,0,0)-tree} with {\bf base} $T_0$, {\bf stages} $T_1, \ldots, T_n$, and {\bf set of connecting colors}  $S_{n}:=\{\delta_0, \delta_1, \ldots, \delta_{n}\}$.
We write $s(T):=n$. Note that $T_{n}$ can be constructed from $G[g]$ in $O(|V(G)|^2k^2)$ time. 
\end{definition}

Let $T_{n+1}$ be an $(n+1,0,0)$-tree as in the above definition. If $T_{n+1}$ is $\phiv$-elementary, we can further grow $T_{n+1}$ to $cl_{\phiv}(T^+_{n+1}(\phiv,\alpha_0,\dd_{n+1}))$ with one more stage as in the above definition. When $T_{n+1}$ is not $\phiv$-elementary, we grow $T_n$ to a tree $T$ in differently such that $V(T) \subseteq V(T_{n+1})$ and $T$ is not $\phiv$-elementary.

\begin{definition}\label{defi-5-sTP0}
Let $n,m$ be non-negative integers. Suppose there exists an $(n+1,0,0)$-tree $T_{n+1}$ with bases $T_0$, stages $T_1, \ldots, T_{n+1}$, and set of connecting colors $S_n=\{\delta_0, \delta_1,\ldots,\dd_{n+1}\}$.
Set $T_{n,0}=T_n$ and,  for each $i\in [m+1]$, let $T_{n,i}=cl_{\phiv}(T_{n, i-1}+g_{n,i-1})$, where $g_{n,i-1}\in \partial T_{n,i-1}\cap \phiv^{-1}(\dd_n)$ is on a $(T_n,\phiv, \alpha_0, \dd_n)$-path $P$  such that a subptah of $P$ from $T_n$ to $g_{n,i-1}$ but not including $g_{n,i-1}$ is contained in $G[T_{n,i-1}]$. We say that $T_{n,m+1}$ is an {\bf $(n, m+1, 0)$-tree} 
with base $T_0$, stages $T_1, \ldots, T_{n}$, {\bf levels} $T_{n,1}, \ldots, T_{n,m+1}$,  and {\bf $\delta_n$-extending edges} $g_{n,0}, g_{n,1}, \ldots, g_{n,m}$. Any prefix $T$ of $T_{n,m+1}$ properly containing $T_{n,m}$ is called an {\bf $(n,m+1,0)$-Tashkinov tree}, and we write $s(T):=n$ and $\ell(T):=m+1$. (In particular, an $(n,m+1,0)$-tree is also an $(n,m+1,0)$-Tashkinov tree.) Observe that $T_{n,m+1}$ can be found in $O(|V(G)|^2k^2)$ time. 
\end{definition}

Since $T_{n+1}$ is not $\phiv$-elementary, there exists an integer $m\ge 0$ such that $T_{n,m}$ is $\phiv$-elementary and $T_{n,m+1}$ is not $\phiv$-elementary.  
We grow $T_{n,m}$ to a tree $T$ differently such that $V(T)\subseteq V(T_{n,m+1})$ and $T$ is not $\phiv$-elementary:
For each connecting color not already missing at previous vertices, we reserve two colors which we try to avoid when growing the tree.

\begin{definition}\label{defi-2-hi}
Let $n,m,q$ be non-negative integers. 
Suppose $T_{n,m+1}$ is an $(n, m+1,0)$-tree with base $T_0$, stages $T_{1},\ldots, T_n$, set of connecting colors
$S_n=\{\dd_0, \dd_1, \ldots, \dd_n\}$, levels $T_{n,1}, \ldots, T_{n,m+1}$,  and  $\delta_n$-extending edges  $g_{n,0}, g_{n,1}, \ldots, g_{n,m}$.
 Let $T_{n,m,0}=T_{n,m}$ and $g_{n,m,0}=g_{n,m}$. For each 
$\dd\in S_n\setminus \phibar(T_{n,m,0})$, let $\Gamma_0(\delta)\subseteq \phibar(T_{n,m,0}) \setminus \phibar(g)$ with $|\Gamma_0(\dd)|=2$, such that $\Gamma_0(\delta)$, $\dd\in S_n\setminus \phibar(T_{n,m,0})$, are pairwise disjoint. 
\begin{itemize}
    \item [(a)] For $i\in [q]$, choose $g_{n,m,i}\in \partial T_{n,m,i}$ with $\phiv(g_{n,m,i})\in \Gamma_{i-1}(\dd^i)$ for some $\dd^i\in S_n\setminus \phibar(T_{n,m,i})$, let $\Gamma_i(\dd)=\Gamma_{i-1}(\dd)$ if $\dd\in S_n\setminus \phibar(T_{n,m,i})$ and $\dd\ne \dd^i$, and let  $\Gamma_i(\dd^i) = (\Gamma_{i-1}(\dd^i)\cup\{\gamma_i\})\backslash\{\phiv(g_{n,m,i})\}$ for some $\gamma_i\in \phibar(T_{n,m,i} -T_{n,m,i-1})\backslash\{\dd_n\}$. 
 \item [(b)] For $i\in [0,q]$, let $T_{n,m, i+1}$ be an ordered tree obtained from $T_{n,m,i}+g_{n,m,i}$ by adding edges, such that\\ {\em (b1)} for each $e\in E(T_{n,m,i+1})\setminus E(T_{n,m,i})$, $\phiv(e)\in \phibar(T_{n,m,i+1}(e)^-)$, and if $\phiv(e)\in \Gamma_i(\dd)$ and $\dd\in S_n\setminus \phibar(T_{n,m,i})$ then $\dd \in \phibar(T_{n,m,i+1}(e)^-)$, and\\  {\em (b2)} $\phiv(\partial T_{n,m,i+1})$ is disjoint from $\phibar(T_{n,m,i+1}) \backslash \bigcup_{\dd\in S_n\backslash\phibar(T_{n,m,i+1})}\Gamma_{i}(\delta)$. 
 \end{itemize}    
We say that 
$T_{n,m,q+1}$ is an {\bf $(n,m,q+1)$-tree} with base $T_0$, stages $T_1, \ldots, T_{n}$, levels $T_{n,1}, \ldots$, $T_{n,m}$, $\delta_n$-extending edges  $g_{n,0}, g_{n,1}, \ldots, g_{n,m}$,  {\bf phases} $T_{n,m,1}, \ldots, T_{n,m,q+1}$, and \\ {\bf phase-extending} edges $g_{n,m,0}, g_{n,m,1},\ldots, g_{n,m,q}$. Write $\Gamma_i=\bigcup_{\dd\in S_n\backslash \phibar(T_{n,m,i})} \Gamma_i(\dd)$ for $i\in [0, q]$. Any prefix $T$ of $T_{n,m,q+1}$ properly containing $T_{n,m,q}$ is called an {\bf $(n,m,q+1)$-Tashkinov tree}. (Thus, an $(n,m,q+1)$-Tree is also an $(n,m,q+1)$-Tashkinov tree.)  Let $s(T)=n$, $\ell(T)=m$, and $p(T)=q+1$. Note that $T_{n,m,q}$ can be found in $O(|V(G)|^2k^2)$ time. 
\end{definition}

{\bf Remark}. When reducing a non--elementary tree to a `smaller' non-elementary tree, we perform Kempe exchange outside $T_{n,m,i}$ for some $i\in [0,q]$ if one of the colors involved is from $\phibar(T_{n,m,q})\cup S_n$.  Moreover, $\Gamma_i$, $i\in [0, q]$, remain the same during the process.


\begin{definition}\label{defi-trunk}
Let $n,m,q$ be non-negative integers, let $T$ be an $(n,m,q+1)$-Tashkiov tree 
with edge ordering $e_1<e_2<\cdots <e_p$, and let $T_{n,m,q}$ be the final phase of $T$.  Let $j$ be minimum such that $e_j\in E(T)\setminus E(T_{n,m,q})$ and $T[\{e_j, \ldots, e_p\}]$ is a path;  then $branch(T):=T[\{e_j, \ldots, e_p\}]$ is called the {\bf branch} of $T$ and its length is denoted by $b(T)$. Let $v$ denote the end of $branch(T)$ that is not the last vertex of $T$; then the {\bf trunk} of $T$, denoted as $trunk(T)$, is obtained from $T$ by deleting $branch(T)-v$ and let $t(T)$ denote the number of edges in $trunk(T)$. 
 The vector $(s(T), \ell(T), p(T), t(T), b(T))$ defines a lexicographic ordering on all Tashkinov trees. 
\end{definition}

We can now give a brief overview of the proof of Theorem~\ref{thm:main},  and we refer the reader to the proof in Section 5 for a more comlete overview. Let $(G,g,\phiv, \alpha_0)$ be a quadruple. 

If $T_0:=cl_\phiv(G[g])$ is not $\phiv$-elementary then we apply Lemma~\ref{lem:recoloring} and obtain a coloring $\pi$ such that $\pi^{-1}([k])=\phiv^{-1}([k])\cup \{g\}$ and $\pi$ is Kempe equivalent to $\phiv$ on $G-g$.  If $\pi^{-1}([k])=E(G)$, stop and output $\pi$ as the desired coloring; otherwise, update $\phiv$ with $\pi$, update $g$ with an edge from $E(G)\setminus \pi^{-1}([k])$, and repeat the process. 

If $T_0$ is $\phiv$-elementary
then we grow an $(n,m,q+1)$-tree $T_{n,m,q+1}$ as in the above definitions, such that $T_{n,m,q}$ is $\phiv$-elementary but $T_{n,m,q+1}$ is not $\phiv$-elementary. Choose the prefix $T$ of $T_{n,m,q+1}$ such that $T$ is not $\phiv$-elementary and, subject to this, $T$ is `minimal' in a lexicographic order. Then $T$ is an $(n,m,q+1)$-Tashkinov tree. We reduce $T$ to Tashkinov tree $T'$ and a coloring $\phiv'$ Kempe equivalent to $\phiv$ (see  Lemma~\ref{lem-final}), such that $T'$ is not $\phiv'$-elementary and 
$(s(T'),\ell(T'),p(T'),t(T'), b(T'))$ is smaller than $ (s(T), \ell(T), p(T), t(T), b(T))$. Update $T$ with $T'$ and $\phiv$ with $\phiv'$, and repeat the process.

\section{Structure-preserving Kempe exchanges}

Throughout this section, let $(G, g,\phiv,\alpha_0)$ denote a quadruple, with $\phiv(G)\subseteq [k]$ where $k:=\max\{\Delta(G)+1, \Gamma(G)\}$. We consider Kempe exchanges that preserve those structure of $(n,m,q+1)$-Tashkinov trees determined by exit and nonexit edges.

For a tree $T$ in $G$ and $\alpha, \beta\in [k]$,  a {\bf $(T, \phiv,\alpha, \beta)$-exit path} is a path $P$ in $G[\phiv^{-1}(\{\alpha,\beta\})]$ with ends $u,v$ such that $u\in V(T)$, $V(P-u)\cap V(T)=\emptyset$, and $\phibar(v)\cap \{\alpha, \beta\}\ne \emptyset$.
The edge of $P$ incident with $u$ is called a {\bf $(T,\phiv,\alpha, \beta)$-exit edge}. Any edge of $\partial T\cap \phiv^{-1}(\{\alpha,\beta\})$ that is not a $(T, \phiv,\alpha, \beta)$-exit edge is called a {\bf $(T,\phiv,\alpha, \beta)$-nonexit edge}. Notice that a $(T, \phiv, \alpha, \beta)$-nonexit edge is on a $(T,\phiv,\alpha,\beta)$-path in $G$.

To grow/reduce trees using Kempe exchanges, we need to use $(T,\phiv, \alpha_0,\delta_n)$-paths. 
To guarantee the existence of such paths, we will show that for any Tashkinov tree $T$ in $G$ there is at most one $(T,\phiv,\alpha_0,\delta_n)$-exit path. 
First, we consider $(n,0,0)$-trees.

\begin{lemma}\label{thm-2-Tn-extension}
Let $n$ be a non-negative integer and $T:=T_n$ be an $(n,0,0)$-tree in $G$ with base $T_0$, stages $T_1, \ldots, T_n$, and set of connecting colors $S_n:=\{\dd_0, \ldots, \dd_{n}\}$.  Let  $\alpha\in \phibar(T)$ and $\delta\in [k]$.
In $O(|V(G)|^3|E(G)|^2)$ time, 
\begin{itemize}
    \item [(i)] one concludes that there is at most one $(T,\phiv, \alpha, \delta)$-exit path, or 
    \item [(ii)] one finds an $(n',0,0)$-tree $T'$ and a coloring $\phiv'$ Kempe equivalent to $\phiv$ such that $n'\le n$ and $T'$ is not $\phiv'$-elementary.
    \end{itemize}
\end{lemma}

\pf If $T$ is not $\phiv$-elementary then (ii) holds with $T':=T$ and $\phiv':=\phiv$. So assume that $T$ is $\phiv$-elementary.  If $\delta\in\phibar(T)$ then $\alpha, \dd\notin \phiv(\partial T)$ since $\alpha\in \phibar(T)$ and $T$ is $\phiv$-closed; 
so there is no $(T,\phiv,\alpha,\dd)$-exit path and we have (i).  Now assume $\delta\notin \phibar(T)$ and let $P_x,P_y$ be distinct $(T,\phiv,\alpha, \delta)$-exit paths, with $V(P_x)\cap V(T)=\{x\}$, $V(P_y)\cap V(T)=\{y\}$, and $y \in V(T(x))$. Assume $x\in V(T_{j} -T_{j-1})$  for some $j\in [0,n]$ with $T_{-1}=\emptyset$. 
The next four paragraphs describe a process for reducing $|V(P_x)|+|V(P_y)|$ and forcing (ii).

Let $\beta\in \phibar(x)$. 
Let $\phiv_0 = \phiv/(G-T,\alpha, \beta)$, which is Kempe equivalent to $\phiv$ as $\alpha,\beta\notin \phiv(\partial T)$.
Then, since $T_i^+(\phiv, \alpha_0, \dd_i) \subseteq T_n=T$ for $i\in [0, n-1]$ (when $n\ge 1$),  
$T$ is the same $(n,0,0)$-tree under both $\phiv$ and $\phiv_0$. Moreover,  $P_x$ and $P_y$ are $(T,\phiv_0, \beta,\dd)$-exit paths, and  $P_x$ is a $(\beta,\dd)$-chain under $\phiv_0$.

Let $\pi_0 = \phiv_0/P_x$. Then $T_{j-1}^+(\pi_0,\alpha_0, \dd_{n-1})=T_{j-1}^+(\phiv,\alpha_0, \dd_{n-1})$, where $T_{-1}^+(\pi_0,\alpha_0, \dd_{-1})=\emptyset$. Let $T^*=T(x)$ if $T_{j-1}^+(\pi_0,\alpha_0, \dd_{j-1})\subseteq T(x)$, and $T^*=T_{j-1}^+(\pi_0,\alpha_0, \dd_{j-1})$ if $\subseteq T(x) \not \subseteq T_{j-1}^+(\pi_0,\alpha_0, \dd_{n-1})$. Then $T^0=cl_{\pi_0}(T^*)$. Since $V(P_x)\cap V(T_n) =\{x\}$, 
the $(T_{i},\pi_0, \alpha_0, \dd_{i})$-paths are precisely the $(T_{i},\phiv_0, \alpha_0, \dd_{i})$-paths, for $i\in [0,j-1]$ (when $j\ge 1$).  Hence, under $\pi_0$,  $T^0$ is an $(j,0,0)$-tree with base $T_0$ and stages $T_1, \ldots, T_{j-1},T^0$. Note $T^0$ can be found in $O(|V(G)|^2k^2)$ time. Let $x^*$ be the end of $P_x$ other than $x$, and $y^*$ be the end of $P_y$ other than $y$. Also, let $yy'$ be the edge of $P_y$ incident with $y$. Note $\pi_0(yy')=\dd\in \pibar_0(x)$; so $y'\in V(T^0)$.

Suppose there exists $z\in \{x^*, y^*\}\cap V(T^0)$. Note $\pibar_0(z)\cap \{\beta, \dd\} \ne \emptyset $. If $\dd\in \pibar_0(z)$ then  $T^0$ is not $\pi_0$-elementary since $\dd\in \pibar_0(x)$. Now assume  $\beta\in \pibar_0(z)$. Then $V(P_x\cup P_y)\subseteq V(T^0)$; in particular, $x^*,y^*\in V(T^0)$. If $\delta \in \pibar_0(\{x^*,y^*\})$ then $T^0$ is not $\pi_0$-elementary (since $\dd\in \pibar_0(x)$). If $\delta \notin \phibar_0(\{x^*,y^*\})$ then  $\beta\in \pibar_0(x^*)\cap \pibar(y^*)$, and again $T^0$ is not $\pi_0$-elementary. So in all cases, (ii) holds with $T':=T^0$ and $\phiv':=\pi_0$ which can be found in $O(|V(G)|^2k)$ time.

Now suppose $\{x^*,y^*\}\cap V(T^0)= \emptyset$.  Then $P_x$ and $P_y$ contain $(T^0, \pi_0,\beta, \dd)$-exit paths, say  $P_{x_1}$ and $P_{y_1}$, with $P_{x_1}\subseteq P_x$, $P_{y_1}\subseteq P_y$, $V(P_{x_1})\cap V(T^0)=\{x_1\}$, and  $V(P_{y_1})\cap V(T^0)=\{y_1\}$.  Note that $y_1\ne y$ since $y'\in V(T^0)$.  
Thus,  $|V(P_{x_1})| + |V(P_{y_1})| < |V(P_x)| + |V(P_y)|$. 

Repeat the above process for $T^0, \pi_0, \dd, \beta$ (as $T,\phiv, \alpha, \dd$) in the preceding four paragraphs. With at most $|V(G)|$ iterations, we arrive at an $(j,0,0)$-tree $T'$ with base $T_0$, stages $T_1,\ldots, T_{j-1}$, $T'$ and a coloring $\phiv'$ Kempe equivalent to $\phiv$, such that $T'$ is not $\phiv'$-elementary; so (ii) holds. Since $|V(P_x)|+|V(P_y)|\le |V(G)|$, the entire process takes $O(|V(G)|^3|E(G)|^2)$ time as $k=O(|E(G)|)$. \qed

Next, we prove a technical lemma about nonexit edges of $(n,0,0)$-trees.

\begin{lemma}\label{thm-2-handle1}
Let $n$ be a non-negative integer and $T_n$ an $(n,0,0)$-tree, with base $T_0$, stages $T_1, \ldots, T_n$, and set of connecting colors $S_n=\{\delta_0,\ldots, \delta_{n}\}$. 
Let $\alpha_1,\alpha_2\in [k]$,  $R$ an  $(\alpha_1,\alpha_2)$-chain, and $U$ any tree in $G$ containing $V(T_n)$, such that $U$ is $\phiv$-elementary, $\phibar(g)\cap \phiv(\partial U)=\emptyset$, $\alpha_1\in \phibar(U)\setminus (\{\alpha_0\}\cup \phiv(\partial U))$, $\alpha_2\notin \phiv(\partial U)$ when $\alpha_2\in \phibar(U)$, $|V(R)\cap V(U))|\le 1$, and both ends of $R$ are outside $T_n$. 
In $O(|V(G)|^3|E(G)|^2$ time,
\begin{itemize}
    \item [(i)] one concludes that 
    the $(T_n, \phiv/R, \alpha_0,\dd_n)$-nonexit edges are the $(T_n, \phiv, \alpha_0,\dd_n)$-nonexit edges; or
     \item [(ii)] one finds an $(n', 0,0)$-tree $T'$ and coloring $\phiv'$ Kempe equivalent to $\phiv$ such that $n' \le n$ and $T'$ is not $\phiv'$-elementary.
     \end{itemize}
\end{lemma} 

 \pf 
 First, we may assume 
 
 {\it Claim} 1. $\alpha_1\in \phibar(g)\backslash\{\alpha_0\}$ and $\alpha_2\in \{\alpha_0, \dd_n\}$.
 
 Suppose $\alpha_1 = \dd_n$. Then $\dd_n  \notin \phiv(\partial U)$ by assumption. Thus, since $\alpha_0\notin \phiv(\partial U)$ (as $\alpha_0\in \phibar(g)$),   
 $T_n^+(\phiv,\alpha_0, \dd_n)\subseteq U$;  so   $|V(R)\cap V(T_n^+(\phiv, \alpha_0, \dd_n))| \le 1$ and, hence, $E(R)\cap E(T_n^+(\phiv, \alpha_0, \dd_n)) = \emptyset$. Therefore, 
 the $(T_n, \phiv/R, \alpha_0, \dd_n)$-paths are the $(T_n, \phiv,\alpha_0, \dd_n)$-paths and, thus, (i) holds. 

 Thus, we may assume $\alpha_1\ne \dd_n$. Suppose $\alpha_1\notin \phibar(g)\backslash\{\alpha_0\}$. Let $\alpha\in \phibar(g)\backslash\{\alpha_0\}$. Then $\alpha_1,\alpha\notin \phiv(\partial T)\cup \{\alpha_0,\dd_n\}$; so $\phiv^*:= \phiv/(G-U, \alpha_1, \alpha)$ is Kempe equivalent to $\phiv$,
the $(T_n, \phiv^*, \alpha_0, \dd_n)$-paths are precisely the  $(T_n, \phiv,\alpha_0, \dd_n)$-paths, and $R$ is an $(\alpha,\alpha_2)$-chain under $\phiv^*$. 
When $n\ge 1$, for $i\in [0, n-1]$, $E(R)\cap E(T_i^+(\phiv, \alpha_0, \dd_n)) = \emptyset$ since $T_i^+(\phiv,\alpha_0,\dd_i)\subseteq T_n$. Therefore, under $\phiv^*$, $T_n$ is also an $(n,0,0)$-tree. Moreover, $\alpha,\alpha_1, U, R$ satisfy the conditions of the lemma.  Hence, (i) or (ii) holds for $\phiv,\alpha_1$ if, and only if, (i) or (ii) holds for $\phiv^*,\alpha$. 

Thus, we may assume that $\alpha_1\in \phibar(g)\backslash\{\alpha_0\}$. If $\alpha_2\notin \{\alpha_0, \dd_n\}$ then (i) holds. So we may assume $\alpha_2\in \{\alpha_0,\dd_n\}$. The above process takes $O(|V(G)|)$ time. \qed

 We can now prove that (i) or (ii) holds without involving the tree $U$. The only conditions we will use are:    
 $\alpha_1\in \phibar(g)\backslash\{\alpha_0\}$, $\alpha_2\in \{\alpha_0, \dd_n\}$, $|V(R)\cap V(T_n)| \le 1$, and both ends of $R$ are outside of $T_n$. 
 
 Let $P_x$ be the unique $(T_n,\phiv, \alpha_0, \dd_n)$-exit-path and let $P_w$ be the unique $(T_n,\phiv/R, \alpha_0, \dd_n)$-exit-path; else, (ii) follows from Lemma~\ref{thm-2-Tn-extension}. 
 Let $x\in V(P_x)\cap V(T_n)$ and $w\in V(P_w)\cap V(T_n)$.  If $x=w$ then (i) holds. Assume $x\ne w$. 
 
We may assume $w\in V(T_n(x))$. To see this, we observe the symmetry between $T_n$ and $\phiv$ and $T_n$ under $\pi:=\phiv/R$: Since both ends of $R$ are outside $T_n$, it follows that $\pibar(v) = \phibar(v)$ for every $v\in V(T_n)$ and $\pi(e) = \phiv(e)$ for every edge $e\in E(G[T_n])$.  For each $i\in [0, n-1]$, $T_i^+(\phiv, \alpha_0, \dd_i) \subseteq T_n $ and $|V(R)\cap V(T_n)| \le 1$; so the $(T_i, \phiv/R, \alpha_0, \dd_i)$-paths are precisely the $(T_i, \phiv,\alpha_0, \dd_i)$-paths. Hence, 
$T_n$, under $\pi$, is also an $(n,0,0)$-tree with base $T_0$, stages $T_1, \ldots, T_n$, and set of connecting colors $S_n=\{\dd_0, \ldots, \dd_{n}\}$. 
Moreover,  $\alpha_1\in \pibar(g)\backslash\{\alpha_0\}$, and so the conditions on $\alpha_1,\alpha_2,R$ are satisfied under $\pi$ as well.  

  We may also assume that $w$ is an end of a $(T_n,\phiv, \alpha_0, \dd_n)$-path, say $Q_w$; otherwise, (ii) follows from Lemma~\ref{thm-2-Tn-extension}.  Since $Q_w$ is not a $(T_n, \phiv/R, \alpha_0, \dd_n)$-path,  we have  $E(R)\cap E(Q_w)\ne \emptyset$. 
 Let $j\in [-1, n-1]$ such that $x\in V(T_{j+1} - T_j)$, where $T_{-1}=\emptyset$. Then $w\in V(T_{j+1}(x))$. We will construct a coloring $\phiv'$ and a $(j+1,0,0)$-tree $T'$ containing $V(Q_w\cup R)$ (see the next two paragraphs) and show that (ii) holds for $T',\phiv'$. Note that $T_n$ is $\phiv$-elementary (since $U$ is $\phiv$-elementary). 

Let $T^0=T_{j+1}$, $\xi_{-1}=\phiv$, $x_0=x$, $\gamma^{-1}=\alpha_0$, $\beta^{-1}=\dd_n$, and $P_{x_0}=P_x$ which is a $(T^0,\xi_{-1},\beta^{-1},\gamma^{-1})$-exit path. 
Note that $T^0$ is $\xi_{-1}$-closed. Let $\beta^0\in \xibar_{-1}(x_0)$ and $\gamma^0\in \phibar(T(x_0))\setminus \{\gamma^{-1}, \alpha_1, \alpha_2,\beta^0\}$ (which exists as $T^0=T_n$ is $\phiv$-elementary and $|VT(^0)|\ge |V(T_0)|\ge 3$). 
Suppose for some integer $i\ge 0$, we have defined $T^i, \xi_{i-1},x_i, \gamma^{i-1},\beta^{i}$, and $P_{x_i}$, such that $T^i$ is a $(j+1,0,0)$-tree under $\xi_{i-1}$ and contains $T_j$ as its $j$th stage, $T^i$ is $\xi_{i-1}$-elementary and $\xi_{i-1}$-closed, $P_{x_i}$ is a $(T^i,\xi_{i-1},\beta^{i-1},\gamma^{i-1})$-exit path, $x_i$ is the end of $P_{x_i}$ contained in $T^i$,   
$\beta^{i-1}\notin \xibar_{i-1}(T^i)$, $\beta^i\in \xibar_{i-1}(x_i)$,  and $\gamma^{i-1}\in \xibar_{i-1}(T^i({x_{i-1}})-T^{i-1}(x_{i-2}))$ if $i\ge 2$. We define $T^{i+1}, 
\xi_i,x_{i+1}, \gamma^i,\beta^{i+1}$ and $P_{x_{i+1}}$ as follows. 

Let $\phiv_i= \xi_{i-1}/(G-T^i, \beta^i, \gamma^{i-1})$.
Since $\beta^i, \gamma^{i-1}\in \xibar_{i-1}(T^i)$, it follows that 
$\beta^i, \gamma^{i-1}\notin  \phiv(\partial T^i)$,   $\phiv_i$ is Kempe equivalent to $\xi_{i-1}$, and $T^i$ is a $(j+1, 0,0)$-tree.  Now $P_{x_i}$ is a $(\beta^i, \beta^{i-1})$-chain under $\phiv_i$. Thus, 
$\pi_i := \phiv_i/P_{x_i}$ is Kempe equivalent to $\phiv_i$, $P_{x_i}$ is a $(\beta^i, \beta^{i-1})$-chain under $\pi_i$, and $\beta^{i-1}\in \pibar_i(x_i)$. Moreover, the edge of $P_{x_i}$ incident with $x_i$, denoted $x_iy_i$, is the only edge in $\partial T^i\cap \pi_i^{-1}(\beta^i)$. 
Choose $\gamma^i\in \pibar_i(T^i(x_i) - T^{i-1}(x_{i-1}))$ which exists as $T^i(x_i)$ is $\pi$-elementary (since $T^i$ is $\phiv$-elementary and $\beta^{i-1}\notin \xibar_{i-1}(T^i)=\phibar(T^i)$), and 
let $Q_{x_i}$ be the $(T^i,\pi_i, \gamma^i, \beta^i)$-exit path (so $V(Q_{x_i}\cap T^i)=\{x_i\}$ and $x_iy_i\in E(Q_{x_i})$). Let $\xi_i=\pi_i/(G-(T^i\cup Q_{x_i}),\beta^i,\gamma^i)$, which is Kempe equivalent to $\pi_i$. Note that $\gamma^i\in \xibar_i(T^i(x_i) - T^{i-1}(x_{i-1}))$ and $\beta^{i-1}\in \xibar_i(x_i)$. Note that both $T^i(x_i)$ and $T_j^+(\xi_i,\alpha_0,\dd_j)$ contain $T_j$ as a prefix; so the edges of $T^i(x_i)\cup T_j^+(\xi_i,\alpha_0,\dd_j)$ can be ordered to form an ordered tree $U_j$ such that $U_j$ contains $T_j^+(\xi_i, \alpha_0,\dd_j)$ as a prefix. Let $T^{i+1}=cl_{\xi_i}(U_j)$. Let $x_{i+1}\in V(Q_{x_i})\cap V(T^{i+1})$, such that  $P_{x_{i+1}}$, the subpath of $Q_{x_i}$ between $x_{i+1}$ and the end $x_i^*$ of $Q_{x_i}$ other than $x_i$, is minimal. 
Note that $P_{x_{i+1}}$ is a $(T^{i+1},\xi_i, \beta^i,\gamma^i)$-exit path. 
Let $\beta^{i+1}\in \xibar_i(x_{i+1})$.

We stop the above process when $x_{i+1}=x_i$ or $\beta^{i}\in \xibar_{i}(T^{i+1})$. Note that this process takes $O(|V(G)|^2k^2)$ time. 

Note that if $x_i^*\in V(T^{i+1})$, we may assume $\beta^i \in \xibar_i(x_i^*) \subseteq \xibar_i(T^{i+1})$;  otherwise, $\gamma^i\in \xibar_i(x_i^*)\cap \xibar_i(T^i(x_i))$ which implies that $T^{i+1}$ is not $\xi_i$-elementary and (ii) holds with $T':=T^{i+1}$ and $\phiv':=\xi_{i+1}$.  Since $T^i(x_i)$ properly contains $T^{i-1}(x_{i-1})$ for $i\ge 1$, the above process must stop, and we assume it stops with $i=m$. Thus, for $i\in [m]$, $T^i$ is a $(j+1, 0,0)$-tree under $\xi_{i-1}$.  

Note that $\phiv_{0}(Q_w\cup P_{x_0}) =\{\beta^{-1}, \beta^{0}\}$. Let $Q_w^*$ be the component of $G[\phiv_{0}^{-1}(\{\beta^{-1},\beta^{0}\})]$ containing $Q_w$.  So $V(Q_w^*)\cap V(P_{x_0})=\emptyset$ as $P_{x_0}$ is a $(\beta^{-1}, \beta^{0})$-chain under $\phiv_{0}$. We will show that $T^m$  is not $\xi_{m-1}$-elementary by considering whether $V(Q_w^*\cup R) \subseteq V(T^m)$.
By definition and by Claim 1, $\phiv_{0}(R) \subseteq \{\alpha_1, \beta^{-1},\gamma^{-1}, \beta^0\}$. Also observe the following 

{\it Claim} 2. (1) $\xi_{0}(Q_w^*\cup R\cup P_{x_0})\subseteq \{\alpha_1, \beta^{-1}, \gamma^{-1}, \beta^0, \gamma^0\}$ and $\alpha_1,\beta^{-1},\gamma^{-1}, \gamma^0\in \xibar_0(T^0(x_0))$, and (2) for $i\in [m]$, the colors involved in Kempe exchanges for obtaining $\xi_{i}$ from $\xi_{i-1}$ are $\beta^{i-1}, \beta^{i},\gamma^{i-1}, \gamma^{i}$. \qed

Let $A_1=\{\alpha_1, \beta^{-1},\gamma^{-1}, \gamma^0\}$ and $A_i=A_{i-1}\cup \{\beta^{i-2}, \gamma^{i-1}\}$ for $i\in [m]\setminus \{1\}$.

{\it Claim} 3. For $i\in [m]$, (1) $A_i\subseteq \xibar_{i-1}(T^{i-1}(x_{i-1}))$, and (2) $\xi_{i-1} (Q_w^*\cup R \cup P_x) \subseteq A_{i} \cup \{\beta^{i-1}\}$.

Observe that Claim 3 for $i=1$ follows from (1) of Claim 2.
 Now suppose Claim 3  holds for some $i$, with $1\le i<m$. Since $T^{i-1}(x_{i-1})\subseteq T^i(x_{i})$ and $\beta^{i-1},\gamma^{i}\in \xibar_{i}(T^{i}(x_{i}))$, we have $A_{i+1}=A_i\cup \{\beta^{i-1},\gamma^{i}\}\subseteq \xibar_{i}(T^i(x_{i}))$; so (1) holds.
 
By (2) of Claim 2, we see that  $\xi_i(Q_w^*\cup R\cup P_x)\subseteq \xi_{i-1}(Q_w^*\cup R\cup P_x) \cup \{\beta^{i-1},\gamma^{i-1},\beta^{i}, \gamma^i\} \subseteq (A_i\cup \{\beta^{i-1}\}) \cup  \{\beta^{i-1},\gamma^{i-1},\beta^{i}, \gamma^i\} = A_{i+1}\cup \{\beta^i\}$.  So (2) of Claim 3 also holds. \qed

{\it Claim} 4. Let $w_1,w_2$ be the ends of $R$. Then (ii) holds  or, for $i\in [2]$, $w_i\in V(T^m)$ implies $\beta^{m-1}\in \xibar_{m-1}(w_i)$. 

Let $r\in [m]$ such that $T^{r}$ contains $w_i$, where $i\in [2]$. We may assume that  $T^r$ is  $\xi_{r-1}$-elementary; otherwise, (ii) holds with $T':=T^r$ and $\phiv':=\xi_{r-1}$. 

Since $\xibar_{-1}(w_i)\cap (\xibar_{-1}(T^0(x_0))\cup \{\beta^{-1}\}) \ne \emptyset$, we have $\xibar_{r-1}(w_i)\cap (\xibar_{r-1}(T^0(x_0))\cup A_r\cup\{\beta^{r -1}\}) \ne \emptyset$ by (2) of Claim 2, 
Notice that  $T^0(x_0) \subseteq T^{r-1}(x_{r-1})$ and $A_r \subseteq \xibar_{r-1}(T^{r-1}(x_{r-1}))$ (by (1) of Claim 3). Therefore, $\beta^{r-1}\in \xibar_{r-1}(w_i)\subseteq \xibar_{r-1}(T^r)$ as $T^r$ is $\xi_{r-1}$-elementary. Since the construction process stops at $T^m$, we must have $r=m$;  so $\beta^{m-1}\in \xibar_{m-1}(w_i)$. \qed

{\it Claim} 5. We may assume $\beta^{r-1}\notin \xibar_{r-1}(T^r)$ for all $r\in [m]$.

Since the construction process stops at $T^m$, $\beta^{r-1}\notin \xibar_{r-1}(T^r)$ for all $r\in [m-1]$. Thus, to prove Claim 5, we assume $\beta^{m-1}\in \xibar_{m-1}(T^m)$.  By Claim 3, $\xi_{m-1}(Q_w^*\cup R\cup P_{x_0})\subseteq A_m\cup \{\beta^{m-1}\}\subseteq \xibar_{m-1}(T^{m-1})\cup \{\beta^{m-1}\} \subseteq \xibar_{m-1}(T^m)$. Hence, $V(R)\subseteq V(T^m)$ and, by Claim 4, we have $\beta^{m-1}\in \xibar_{m-1}(w_1)\cap \xibar_{m-1}(w_2)$; so (ii) holds with $T':=T^{m}$ and $\phiv':=\xi_{m-1}$.  \qed

 {\it Claim} 6. $E(Q_w^*\cup R)\cap \partial T^m=\{x_{m}y_{m}\}$. 

Let $uv\in E(Q_w^*\cup R)\cap \partial T^m$ with 
$u\in V(T^m)$ and $v\notin V(T^m)$. By (2) of Claim 3, $\xi_{m-1}(uv) \in A_m\cup \{\beta^{m-1}\}$; so, by (1) of Claim 3, $\xi_{m-1}(uv) =\beta^{m-1}$ since $T^m$ is $\xi_{m-1}$-closed. 
Note that $\pi_{m-1}(uv)\ne \gamma^{m-1}$ since $\gamma^{m-1}$ is not used in any Kempe exchange for obtaining $\pi_{m-1}$  and $\gamma^{m-1}\notin  \{\alpha_1, \alpha_0, \beta^{-1},\gamma^{-1}\}\supseteq \phiv_{0}(Q_w^*\cup R)$.   
Thus, since $\xi_{m-1}=\pi_{m-1}/(G-(T^{m-1}\cup Q_{x_{m-1}}),\beta^{m-1},\gamma^{m-1})$, we have $\xi_{m-1}(uv) = \pi_{m-1}(uv) = \beta^{m-1}$. It follows that $uv \in E(G[T^{m-1}\cup Q_{x_{m-1}}])$. Hence $\{u,v\}\subseteq V(T^{m-1})$ or $uv\in E(Q_{x_{m-1}})$.

Suppose $\{u,v\}\subseteq V(T^{m-1})$. 
Then $\phiv_{m-1}(uv)=\xi_{m-2}(uv)$ since $\phiv_{m-1}=\xi_{m-2}/(G-T^{m-1}, \beta^{m-1},\gamma^{m-2})$. By (2) of Claim 3, $\xi_{m-2} (Q_w^*\cup R) \subseteq A_{m-1} \cup \{\beta^{m-2}\}$. Thus, $\phiv_{m-1}(uv)\in A_{m-1}\cup \{\beta^{m-2}\}$; so $\phiv_{m-1}(uv) \ne \beta^{m-1}$. Since obtaining $\pi_{m-1}$ and $\xi_{m-1}$ from $\phiv_{m-1}$ involve no edges with both ends in $V(T^{m-1})$,  $\xi_{m-1}(uv)\ne \beta^{m-1}$, a contradiction. 

Hence, $uv\in E(Q_{x_{m-1}})$. By Claim 5, $x_m=x_{m-1}$. Thus, $x_m$ is both the first vertex  of $P_{x_{m}}$ and the last vertex of $Q_{x_{m-1}}$ contained in $V(T^m)$; so $uv=x_my_m$. \qed

 We may assume $\{w_1,w_2\} \not \subseteq V(T^m)$; for otherwise, $\beta^{m-1}\in \xibar_{m-1}(w_1)\cap\xibar_{m-1}(w_2)$ (by Claim 4) and, hence,(ii) holds with $T':=T^m$ and $\phiv':=\xi_{m-1}$. Therefore, since $w\in V(T^0(x_0))\subseteq V(T^m)$, 
there exists $uv\in E(Q_w^*\cup R)\cap \partial T^m$. By Claim 6, $uv=x_my_m$ is the unique edge in $E(Q_w^*\cup R)\cap \partial T^m$. Thus,  $V(Q_w^*)\subseteq V(T^m)$ as $Q_w^*$ is a cycle.  Hence, $x_my_m\in E(R) \setminus E(Q_w^*)$ and an end of $R$, say $w_1$, is in  $T^m$. By Claim 4, we may assume $\beta^{m-1}\in \xibar_{m-1}(w_1)$ (or (ii) holds). Thus, $V(Q^*_w\cup R) \subseteq V(T^m)$ (by Claim 3), forcing $\{w_1,w_2\}\subseteq  V(T^m)$, a contradiction. 

The complexity of the proof is $O(|V(G)|^3|E(G)|^2)$ time, due to the use of Lemma~\ref{thm-2-Tn-extension}. \qed

Next we show that $(n,m,q)$-trees have at most one exit path.
\begin{lemma}\label{thm-2-T(n,m)-exit} 
Let $n,m,q$ be non-negative integers, $T$ be an $(n,m,q)$-tree, and $\alpha, \delta\in [k]$, such that $\alpha\in \phibar(T)\backslash\{\alpha_0\}$ and $\alpha\notin \phiv(\partial T)$. 
In $O(|V(G)|^4|E(G)|^2)$ time,  
\begin{itemize} 
  \item [(i)] one concludes that there is at most one $(T,\phiv, \alpha, \delta)$-exit path; or 
  \item [(ii)] one finds an $(n',m',q')$-tree $T'$ and a coloring $\phiv'$ Kempe equivalent to $\phiv$ such that $(n',m',q')\le (n,m,q)$ and $T'$ is not $\phiv'$-elementary. 
  \end{itemize} 
\end{lemma}

\pf We may assume that $T$ is $\phiv$-elementary; for otherwise, (ii) holds with $T':=T$ and $\phiv':=\phiv$.  By Lemma~\ref{thm-2-Tn-extension}, we may assume $(m,q) \ne (0,0)$.  

Let $T_0$ be the base of $T$, let $T_1, \ldots, T_n$ be the stages of $T$, let $T_{n,1},\ldots, T_{n,m}$ be the levels of $T$, and let $T_{n,m,1},\ldots, T_{n,m,q}$ be the phases of $T$. Then $T=T_{m,n,q}$. Let $S_n=\{\dd_0, \dd_1,\ldots,\dd_n\}$ denote the set of connecting colors of $T$. 

We may assume $\alpha\in \phibar(g)\backslash\{\alpha_0\}$.  For, otherwise, choose $\alpha^*\in \phibar(g)\backslash\{\alpha_0\}$. Since $\alpha^*, \alpha\notin \phiv(\partial T)$, $\phiv^* := \phiv/(G-T, \alpha, \alpha^*)$ is Kempe equivalent to $\phiv$. Clearly,  $T_i^+(\phiv^*, \alpha_0, \dd_i) = T_i^+(\phiv, \alpha_0, \dd_i)$ for $i\in [0, n-1]$ (as they are inside $T$). Note that if $\alpha=\dd_n$ then $V(T_n^+(\phiv, \alpha_0, \dd_n))\subseteq V(T)$ as $\alpha,\alpha^*\notin \phiv(\partial T)$, and if $\alpha \ne \dd_n$ then $\alpha, \alpha^*\notin \{\alpha_0,\dd_n\}$.  
Hence, the $(T_n,\phiv^*, \alpha_0, \dd_n)$-paths are precisely the $(T_n,\phiv, \alpha_0, \dd_n)$-paths. 
Thus, under $\phiv^*$, $T$ is also an $(n,m,q)$-tree. 
 Moreover, there is at most one $(T,\phiv^*, \alpha^*, \delta)$-exit path if, and only if, there is at most one $(T,\phiv, \alpha, \delta)$-exit path.
 So could  work with $\phiv^*,\alpha^*$ instead of $\phiv,\alpha$.

To prove the lemma, suppose there exist two $(T,\phiv, \alpha, \delta)$-exit paths, say $P_x$ and $P_y$, with $V(P_x)\cap V(T)=\{x\}$ and $V(P_y)\cap V(T)=\{y\}$. Without loss of generality, let $y\in V(T(x))$. We may assume $x\notin V(T_n)$ by Lemma~\ref{thm-2-Tn-extension}. Let $(m',q')$ be maximal with respect to the lexicographic ordering subject to $x\notin V(T_{n,m',q'})$. We will prove (ii) by repeating the process in the next five paragraphs. 


Let $\beta\in \phibar(x)$. Then, since $T$ is $\phiv$-elementary,  
$\beta\notin \Gamma_{q'}\cup\{\alpha_0\}$; so $\alpha, \beta\notin \phiv(\partial T_{n,m,q'})$. Hence, $\phiv_0 := \phiv/(G-T, \alpha, \beta)$ is Kempe equivalent to $\phiv$. Moreover, for $i\in [0, n-1]$, $T_i^+(\phiv, \alpha_0, \dd_i)\subseteq T_n \subseteq T(x)^-$; so $T_i^+(\phiv_0, \alpha_0, \dd_i) = T_i^+(\phiv, \alpha_0, \dd_i)$. Note that $\beta=\dd_n$ implies $T_n^+(\phiv, \alpha_0, \dd_n)\subseteq T_{n,m',q'}$ and that $\beta\ne \dd_n$ implies $\alpha, \beta\notin \{\alpha_0, \dd_n\}$; thus $T_n^+(\phiv_0, \alpha_0, \dd_n) = T_n^+(\phiv, \alpha_0, \dd_n)$. Hence, $T(x)$ has the same tree structure under $\phiv$ and $\phiv_0$.  

Now  $P_x$ is a $(\beta,\dd)$-chain under $\phiv_0$ and $P_y$ is a $(T, \phiv_0, \beta, \dd)$-exit path. Hence, $\pi_0: = \phiv_0/P_x$ is Kempe equivalent to $\phiv_0$, $\dd\in \pibar_0(x)$, and $P_x$ is a $(\beta,\dd)$-chain under $\pi_0$. Since $x\notin V(T_{n,m',q'})$,  both ends of $P_x$ are outside $V(T_n)$. 
We may assume that 
the $(T_n,\pi_0, \alpha_0, \dd_n)$-nonexit edges are precisely the $(T_n,\phiv_0,\alpha_0, \dd_n)$-nonexit edges.  
For, otherwise, by applying Lemma~\ref{thm-2-handle1} to  $T_n,\phiv_0,\alpha, \beta,P_x$, there exist an $(n',0,0)$-tree $T'$ and a coloring $\phiv'$ Kempe equivalent to $\phiv_0$ (and hence $\phiv$), such that $n'\le n$ and $T'$ is not $\phiv'$-elementary and (ii) holds.

Let $T^1$ be an ordered tree obtained from $T(x)$ by successively adding edges $e$ such that (1) $\pi_0(e)\in \pibar_0(T^1(e)^-)$, and if $\pi_0(e)\in \Gamma_{q'}(\dd')$ and $\dd'\in S_n\setminus \pibar_0(T_{n,m',q'})$ then $\dd'\in \pibar_0(T^1(e)^--T(x))$, and (2) $\pi(\partial T^1)$ is disjoint from $\pibar(T^1)$ and $\bigcup_{\dd'\in S_n\setminus \pibar_0(T_{n,m',q'})}\Gamma(\dd')$. Thus,$T_{n,m',q'} \subseteq T(x) \subseteq T^1$. 
It follows that  $T^1$ is an $(n,m',q'+1)$-tree under $\pi_0$ and has the same base, stages, and levels as $T_{n,m',q'}$ under $\phiv$ but with one more phase. Note that $T^1$ contains the edge of $P_y$ incident with $y$.  Let $x^*$ and $y^*$ be the ends of $P_x$ and $P_y$ in $V(T^1)$, respectively. Recall that $\dd\in \pibar_0(x)$, $\pibar_0(x^*)\cap \{\beta,\dd\}\ne \emptyset$, $\pibar_0(y^*)\cap \{\beta,\dd\}\ne \emptyset$, and $\pi_0(P_x\cup P_y)=\{\beta,\dd\}$.  

Suppose there exists $z\in \{x^*,y^*\}\cap V(T^1)$. Note that $\{\beta,\delta\}\cap  \pibar_0(z)\ne \emptyset$. If $\delta\in \pibar_0(z)$ then $T^1$ is not $\pi_0$-elementary. Suppose $\beta\in \pibar_0(z)$; then $\{x^*,y^*\}\subseteq V(T^1)$. If $\delta\in \pibar_0(\{x^*,y^*\})$ then  $T^1$ is not $\pi_0$-elementary; otherwise $\beta\in \pibar_0(x^*)\cap \pibar_0(y^*)$ and again  $T^1$ is not $\pi_0$-elementary. Thus, (ii) holds with $T':=T^1$ and $\phiv':=\pi_0$.

Now assume that $\{x^*,y^*\}\cap V(T^1)= \emptyset$.
Then, $P_x$ and $P_y$ each contain a $(T^1,\pi_0, \beta, \dd)$-exit path. Denote these two subpaths by $P_{x_1}$ and $P_{y_1}$, where $x_1\in V(P_{x_1})\cap V(T^1)$ and $y_1\in V(P_{y_1})\cap V(T^1)$. Since $y\in V(T^1)$, 
 $|V(P_{x_1})| + |V(P_{y_1})| < |V(P_x)| + |V(P_y)|$. 

 Note that the process in the above 5 paragraphs takes $O(|V(G)|^3|E(G)|^2)$ time (due to an application of Lemma~\ref{thm-2-handle1}) and it is repeated at most $|V(G)|$ times to arrive at (ii). Hence, the complexity of the proof is $O(|V(G)|^4|E(G)|^2$.  \qed

We will make heavy use of Kempe exchanges satisfying the conditions of the lemma below. 

\begin{lemma}\label{lem-7-Pt-disjoint}
Let $n,m,q$ be non-negative integers and $T$ an $(n,m,q+1)$-Tashkinov tree, with set of connecting colors $S_n=\{\dd_0, \ldots, \dd_n\}$, and final phase $T_{n,m,q}$. 
Let $v, w\in V(T)$ with $v\in V(T(w))$, and let $\epv_1\in (\phibar(v)\setminus \{\alpha_0\})\setminus \phiv(T(w) - E(T_{n, m, q}))$ and $\epv_2\in \phibar(w)$ such that $\{\epv_1, \epv_2\}\cap \phibar(T_{n,m, q})\not\subseteq \{\alpha_0\}$ or $\{\epv_1, \epv_2\}\cap (\phibar(T_{n,m,q})\cup S_n)=\emptyset$. 
In $O(|V(G)|^4|E(G)|^2)$ time, 
\begin{itemize}
    \item [(i)] one concludes that, for any $(\epv_1, \epv_2)$-chain $R$ not containing $v$, \\  {\em (i1)} 
    the  $(T_n,\phiv/R, \alpha_0, \dd_n)$-nonexit edges are precisely the  $(T_n,\phiv, \alpha_0, \dd_n)$-nonexit edges,\\ 
    {\em (i2)}  $T_{n,m,q}$ is an $(n,m,q)$-tree under $\phiv/R$, \\ 
    {\em (i3)} if $w\notin V(T_{n,m,q})$ then $T(w)$ is an $(n,m,q+1)$-Tashkinov tree under $\phiv/R$, \\
    {\em (i4)} 
    $v,w$ are the ends of an $(\epv_1, \epv_2)$-chain under $\phiv$, and\\
    {\em (i5)} if $w\notin V(T_{n,m,q}\cup R)$ and if $\epv_1\notin \Gamma_q$, or $\epv_1\in \Gamma_q$ and $\epv_1, \epv_2\notin \phiv(T- E(T(w)))$,  or $\epv_1\in \Gamma_q(\dd)$ for some $\dd \in \phibar(T(w))$, then $T$ is an $(n, m, q+1)$-Tashkinov tree under $\phiv/R$; or    
    
    \item [(ii)] one finds an $(n',m',q')$-Tashkinov tree $T'$ and  a coloring $\phiv'$ Kempe equivalent to $\phiv$, such that $T'$ is not $\phiv'$-elementary, and $T'\subseteq T(w)$  or $(n', m', q') < (n,m, q+1)$. 
    \end{itemize}
\end{lemma}

\proof 
We remark that the complexity $O(|V(G)|^4|E(G)|^2)$ is due to the use of 
Lemma~\ref{thm-2-T(n,m)-exit}.  
We may assume that $T^-$ is $\phiv$-elementary; else (ii) holds with $\phiv':=\phiv$ and $T':=T^-$. Let $R$ be an arbitrary  $(\epv_1, \epv_2)$-chain not containing $v$.

We may further assume $v\in V(T_{n,m,q})$. For, otherwise,
we have $\{\epv_1, \epv_2 \}\cap (\phibar(T_{n,m,q})\cup S_n) =\emptyset$ (by assumption) and, hence, (i1) and (i2) hold.  
Since $\epv_1\notin \phiv(T(w) - E(T_{n,m,q}))$ and $\epv_2\in \phibar(w)$, 
we have $E(R)\cap E(T(w)) = \emptyset$; so (i3) holds. Note that $\epv_1,\epv_2\notin \Gamma_q$ as $v, w\notin V(T_{n,m,q})$ and that $w\notin V(R)$ implies $\pibar(w) = \phibar(w)$. Hence, if $w\notin V(R)$ then $T$ is a prefix of a closure of $T(w)$ under $\phiv/R$ (see (b1) of Definition~\ref{defi-2-hi}); thus (i5) holds. 
Now suppose (i4) fails. Then $R'$, the $(\alpha,\beta)$-chain containing $w$, does not contain $v$. The same argument above for $R$ also applies to $R'$; hence (i3) holds for $\xi:=\phiv/P_w(\phiv, \alpha, \beta)$. 
Now $\epv_1\in \xibar(v)\cap \xibar(w)$; so (ii) holds with $T'=T(w)$ and $\phiv':=\xi$.

Thus, $\epv_1\in \phibar(T_{n,m, q})\backslash\{\alpha_0\}$. Note that for any $\gamma\in \{\epv_1,\epv_2\}\setminus \{\alpha_0\}$ with $\gamma\in \phibar(T_{n,m,q})$, $\gamma\notin \phiv(\partial T_{n,m,h})$ for some $h\in [0, q]$; since if $\gamma\in \phibar(T_{n,m,0})$ then $\gamma\notin \phiv(\partial(T_{n,m,0})$, and if $\gamma\in \phibar(T_{n,m,h} -T_{n,m,h-1})$ for some $h \ge 1$ then $\gamma\notin \phiv(\partial(T_{n,m,h}))$. Choose such a $\gamma$ to maximize $h$. We may assume that there is at most one $(T_{n,m,h},\phiv, \epv_1, \epv_2)$-exit path (else (ii) follows from Lemma~\ref{thm-2-T(n,m)-exit}), and if there is one then it must be the $(\epv_1,\epv_2)$-chain $P_v$ under $\phiv$ containing $v$. Therefore,  $R\cap T_{n,m,h} = \emptyset$ as $v\notin V(R)$.  So it follows from the maximality of $h$ that, for each $i\in [q]\setminus [h]$, $\{\epv_1,\epv_2\}\subseteq \phiv(\partial T_{n,m,i})$ and, thus, $\epv_1\in \Gamma_{i-1}$, and  if $\epv_2\in \phibar(T_{n,m,i-1})$ then  $\epv_2\in \Gamma_{i-1}$. 

Let $\eta\in \{\epv_1, \epv_2\}\setminus \{\gamma\}$.
If $\eta\notin \phibar(T_{n,m,h})$ or $\eta\in \phibar(T_{n,m,h})\backslash\phiv(\partial T_{n,m,h})$ then, applying Lemma~\ref{thm-2-handle1} to $T_n$ and $R$ (with $U=T_{n,m,h}$), we see that (i1) or (ii) holds. 
Now assume that $\eta\in \phibar(T_{n,m,h})\cap  \phiv(\partial (T_{n,m,h}))$. Then $\eta\in \Gamma_{h-1}(\dd)$ for some $\dd\in S_n \setminus \phibar(T_{n,m,h})$; hence, $\eta\notin \{\alpha_0, \dd_n\}$.  If $\gamma\ne \dd_n$ then $\{\gamma, \eta\}\cap \{\alpha_0, \dd_n\} = \emptyset$; so (i1) holds. If $\gamma=\dd_n$ then $\dd_n\notin \partial T_{n,m,h}$ which implies $V(T_n^+(\phiv,\alpha_0,\dd_n)) \subseteq V(T_{n,m,h})$; hence, $V(R)\cap V(T^+_n)\subseteq V(R)\cap V(T_{n,m,h}) = \emptyset$ and, clearly, (i1) holds.

Therefore, to complete the proof of this lemma, we suppose (i1) holds and show that (i2)--(i5) hold or (ii) holds.  For convenience, let $\pi=\phiv/R$. Since $R\cap T_{n,m,h} = \emptyset$ and (i1) holds, $T_{n,m,h}$ is an $(n,m, h)$-tree under $\pi$. For each $i\in [q]\setminus [h]$, we have $\gamma\in \phiv(\partial T_{n,m,i})$ and, hence, $\gamma\in \Gamma_{i-1}$. 

{\it Case 1}. $w\in V(T_{n,m,q})$.   

Then, for each $i\in [q]\setminus [h]$, 
$\{\epv_1, \epv_2\}\subseteq \phiv(\partial T_{n,m,i})$; so $\epv_1, \epv_2\notin \phiv(T_{n,m,i} - E(T_{n,m,i-1}))$.
Moreover, $\epv_1, \epv_2\notin \phibar(T_{n,m,i})\backslash\bigcup_{\dd\in S_n\backslash\phibar(T_{n,m,i})}\Gamma_{i-1}(\dd)$.
Hence, if $w\notin V(R)$ then, for $i\in [q]\setminus [h]$, $T_{n,m,i}$ is an $(n, m, i)$-tree under $\pi$ and, hence, (i2) holds; and if $w\in V(R)$ then $T(w)$ is an $(n,m, j+1)$-Tashkinov tree and $\epv_1\in \pibar(v)\cap \pibar(w)$ where  
$w\in V(T_{n,m,j+1} -T_{n,m,j})$ for some $j\in [q-1]$  and, thus, (ii) holds.
Notice that (i3) and (i5) do not apply in this case. For (i4), suppose $w\in V(T^-)$ and $v,w$ are not the ends of an $(\epv_1,\epv_2)$-chain. If $w\in V(T_{n,m,h})$ then there are two $(T_{n,m,h}, \phiv, \epv_1, \epv_2)$-exit paths; hence, (ii) holds by Lemma~\ref{thm-2-T(n,m)-exit}. If $w\notin V(T_{n,m,h})$, say $w\in V(T_{n,m,j+1} -T_{n,m,j})$ for some $j \ge h$. Since $T_{j+1}$ is an $(n, m,j+1)$-tree under $\pi$ (by (i2)), $T(w)$ is an $(n,m,j+1)$-Tashkinov tree under $\pi$. Now $\epv_1\in \pibar(v)\cap \pibar(w)$;
so  (ii) holds with $T':=T(w)$ and $\phiv':=\pi$.  

{\it Case} 2. $w\notin V(T_{n,m,q})$. 

Then $\gamma=\epv_1$ and $\eta=\epv_2$. 
For each $i\in [q]\setminus [h]$, 
we have $\epv_1\notin \phibar(T_{n,m,i})\backslash \Gamma_{i-1}$ as $\epv_1\in \Gamma_{i-1}$. Hence,   $\epv_1\notin \phiv(T(w) - E(T_{n,m,h}))$.  Since $\epv_2\in \phibar(w)$, it follows that $\epv_2\notin \phiv(T(w) - E(T_{n,m,h}))$, unless $h=0$ and $\epv_2=\dd_n$ in which case $\phiv^{-1}(\epv_2)\cap (E(T(w))\setminus E(T_{n,m,h}))) =\{g_{n,m,0}\}$. Hence, since $R\cap T_{n,m,h} = \emptyset$, we have $E(R)\cap E(T(w)- E(T_{n,m,h})) = \emptyset$. Therefore, $T(w)$ satisfies (b1) of Definition~\ref{defi-2-hi} under $\pi$. Moreover, since $\epv_1, \epv_2 \notin \phibar(T_{n,m,i})\backslash \Gamma_{i-1}$ for $i\in [q]\setminus [h]$, $T_{n,m,i}$ under $\pi$ satisfies (b2) of Definition~\ref{defi-2-hi}. Thus, $T(w)$ is an $(n,m, q+1)$-Tashkinov tree. Note that $\phibar(v) = \pibar(v)$ for every $v\in V(T(w))\backslash\{w\}$.  So (i3) holds, which also implies (i2). 

Now suppose that the conditions of (i5) are satisfied.  Then $\pibar(T(w)) = \phibar(T(w))$ since $w\notin V(R)$;  so $T$ is a prefix of some closure of $T(w)$ under $\pi$. Note $\epv_2\notin \Gamma_q$ as $\epv_2\notin \phibar(T_{n,m,q})$. To prove (i5), it suffices to verify that $T$ and $\pi$ also satisfy (b1) of Definition~\ref{defi-2-hi}. But this can be easily checked since $\epv_1\notin \Gamma_q$, or $\epv_1, \epv_2\notin \phiv(T- E(T(w)))$ (hence $\epv_1, \epv_2\notin \pi(T- E(T(w)))$), or $\epv_1\notin \Gamma_q(\dd)$ for some $\dd\in \phibar(T(w))$.
For (i4), suppose that $v,w$ are not the ends of an $(\epv_1,\epv_2)$-chain under $\phiv$. By the same argument as in the end of Case 1, we see that (ii) holds with $T':=T(w)$ and $\phiv':=\pi$. 
\qed

{\bf Remark}. We note that when $w\in V(T^-)$ the tree $T'$  in (ii) of Lemma~\ref{lem-7-Pt-disjoint} satisfies\\ $(s(T'),\ell(T'),p(T'),t(T'), b(T')) < (s(T),\ell(T), p(T), t(T), b(T))$. 

We need to use special colors for Kempe exchanges in order to preserve base, stages, levels and phases. 
To guarantee such colors are available, we use a result of Stiebitz et al.~\cite{StiebitzBook} which can be phrased as: $|V(T_0)|\ge 7$ or there exists a coloring $\phiv'$ Kempe equivalent to $\phiv$ such that $G[g]$ is not $\phiv'$-elementary.

\begin{lemma}\label{lem-7-B} 
Let $n,m,q$ be non-negative integers and let  $T$ be an $(n,m,q+1)$-Tashkinov tree with final phase $T_{n,m,q}$ and set of connecting colors $S_n$. Let $u_r$ be the last vertex of $trunk(T)$, and let $B(T(u_r)^-, \phiv)$ consist of the colors in $\phibar(T(u_r)^-)$ but not in any of the following sets: $\phibar(g)$, $\phiv(T(u_r) - E(T_n))$, $\cup_{\dd\in S_n \setminus \phibar(T(u_r)^-)}\Gamma_q(\dd)$, $S_n \setminus \phibar(T_{n,m,q})$. Suppose $T(u_r)$ is $\phiv$-elementary and $T_{n,m,q}\subseteq T(u_r)$.
Then  $|B(T(u_r)^-, \phiv)| \ge 2$. 
\end{lemma}  

\pf  Since $T(u_r)^-$ is $\phiv$-elementary,  we have
$$|\phibar(T(u_r)^-)| = |\phibar(g)|+|\phibar(T_0)\setminus \phibar(g)| + |\phibar(T_n - T_0)| + |\phibar(T(u_r)^- -T_n)|.$$ 
Note that $|\phibar(v)| \ge 1 $ for each $v\in V(G)$.  Hence, 
\begin{itemize}

\item $|\phibar(T_0)\setminus \phibar(g)|\ge |V(T_0)| -2 \ge 5$ (since $|V(T_0)| \ge 7$); 

    \item $|\phibar(T_n-T_0)| \ge |V(T_n -T_0)| \ge 2n$ (since $|V(T_{i+1} -T_i)| \ge 2$ for $i\in [0,n-1]$); 
    \item $|\phibar(T(u_r)^- -T_n)| \ge  |V(T(u_r) -T_n)| -1 \ge |\phiv(T(u_r)-E(T_n))|-1$;
    \item $|\cup_{\dd\in S_n\setminus \phibar(T(u_r)^-)}\Gamma_q(\dd)| \le 2|S_n \setminus \phibar(T(u_r)^-)|$ and \\$2|S_n \setminus  \phibar(T(u_r)^-)| +  |(S_n\backslash \phibar(T_{n,m,q}))\cap \phibar(T(u_r)^-)| \le 2|S_n| \le 2(n+1)$.
    \end{itemize}
By a simple calculation, we have $|B(T(u_r)^-, \phiv)| \ge 2$.
    \qed

{\bf Remark}. For later convenience, we mention two properties of $\beta\in B(T(u_r)^-, \phiv)$: if $\beta\notin \phibar(T_{n,m,q})$ then $\beta\notin S_n$; if $\beta\in \Gamma_q(\dd)$ then $\dd\in \phibar(T(u_r)^- - T_{n,m,q})$. 

\section{Reducing non-elementary trees}

In this section, we reduce a non-elementary tree a smaller non-elementary tree. Let $(G,g,\phiv,\alpha_0)$ be a quadruple. For convenience, we let $P_v(\alpha,\beta;\phiv)$ to denote the $(\alpha,\beta)$-chain in $G$ containing the vertex $v$ under the coloring $\phiv$.

Let $T$ be an $(n,m,q+1)$-Tashkinov tree, such that $T$ is not $\phiv$-elementary, but $T^-$ is $\phiv$-elementary. For each color $\epv\in \phibar(T^-)$, denote by $z_\epv$ the unique vertex in $V(T^-)$ with $\epv\in \phibar(z_\epv)$.  Let $V(branch(T))=\{v_0,v_1,\ldots v_t\}$, let $f_i$ be the edge of $branch(T)$ between $v_{i-1}$ and $v_i$ for $i\in [t]$, and let $\alpha\in \phibar(v_t)\cap \phibar(T^-)$. 

Let $T_0$ be the base of $T$, $\, T_1, \ldots, T_n$ be the stages of $T$, $\, T_{n,1}, \ldots, T_{n,m}$ be the levels of $T$, and $T_{n,m,1}, \ldots, T_{n,m,q}$ be the phases of $T$. Let $S_n=\{\dd_0, \ldots, \dd_n\}$ be the set of connecting colors of $T$, $g_{n,0}, g_{n,1},\ldots, g_{n,m}$ be the $\dd_n$-extending edges, and $g_{n,m,0}, g_{n,m,1},\ldots, g_{n,m,q}$ denote the phase-extending edges.  Recall $\Gamma_i$, $i\in [q]$, from Definition 2.3.

\begin{lemma}\label{lem-7-1st-edge-elem} 
In $O(|V(G)|^4|E(G)|^2)$ time, (i) one concludes that $T\ne T_{n,m,q} +g_{n,m,q}$, or (ii) one finds a tree $T'$ and a coloring $\phiv'$ Kempe equivalent to $\phiv$, such that 
$T'$ is not $\phiv'$-elementary and $(s(T'),\ell(T'),p(T'),t(T'), b(T')) <  (s(T),\ell(T),p(T),t(T), b(T))$.
\end{lemma}

\proof Suppose $T=T_{n,m,q} +g_{n,m,q}$; so $t=1$, $f_1=g_{n,m,q}$, and $v_1\notin V(T_{n,m,q})$.  
 
  We may assume that  $\alpha\in  \phibar(g)\backslash\{\alpha_0\}$. Otherwise, suppose $\alpha=\alpha_0$ or $\alpha\notin \phibar(g)$.
 Let $\alpha'\in \phibar(g)\backslash\{\alpha_0\}$ and $R:=P_{v_1}(\alpha, \alpha';\phiv)$. By definition, $\alpha'\notin \Gamma_{q-1}$, and so  $\alpha'\notin \partial(T_{n,m,q})$.  
 We may assume that there is at most one $(T_{n,m,q},\phiv,\alpha, \alpha')$-exit path (or (ii) holds by Lemma~\ref{thm-2-T(n,m)-exit}); so $R\cap T_{n,m,q}=\emptyset$. Apply Lemma~\ref{lem-7-Pt-disjoint} with $\{\epv_1,\epv_2\}=\{\alpha',\alpha\}$, $\{v,w\}=\{z_{\alpha'},z_{\alpha}\}$, and $v\in V(T(w))$.  We may assume that   $(T_n, \phiv/R, \alpha_0,\dd_n)$-nonexit edges are precisely the $(T_n, \phiv, \alpha_0,\dd_n)$-nonexit edges, and $T_{n,m,q}$ is the same $(n,m,q)$-tree under both $\phiv/R$ and $\phiv$;  otherwise (ii) holds by (i1) and (i2) of Lemma~\ref{lem-7-Pt-disjoint}. 
Hence, we could consider $\phiv/R,\alpha'$ instead of $\phiv,\alpha$.

Suppose $q\ge 1$. Then $\gamma:=\phiv(g_{n,m,q})\in \Gamma_{q-1}$. Let $R:=P_{v_t}(\alpha, \gamma;\phiv)$.  Apply Lemma~\ref{lem-7-Pt-disjoint} with $\{\epv_1,\epv_2\} =\{ \alpha,\gamma\}$, $\{v,w\}=\{z_{\alpha},z_{\gamma}\}$, and $v\in V(T(w))$. We may assume that $z_\alpha, z_\gamma\notin V(R)$ (or (ii) holds by (i4) of Lemma~\ref{lem-7-Pt-disjoint}).   Hence, $R$ contains two $(T_{n,m,q}, \phiv,\alpha, \gamma)$-exit paths; so (ii)  follows from Lemma~\ref{thm-2-T(n,m)-exit}.

Now assume $q=0$. Then $g_{n,m,q}=g_{n,m}$ is a $\dd_n$-extending edge, and some $(T_n,\phiv,\alpha_0, \dd_n)$-path has a subpath, say $Q_1$, starting from $T_n$ and ending with the edge $g_{n,m,q}=v_0v_1$ such that $V(Q_1-v_1) \subseteq V(T_{n,m})$.    Let $R:=P_{v_1}(\alpha, \alpha_0;\phiv)$. Note that $R\cap T_{n,m,q}=\emptyset$ since $\alpha, \alpha_0\notin \phiv(\partial T_{n,m,q})$.  
Apply Lemma~\ref{thm-2-handle1} to $T_{n}, \phiv, \alpha,\alpha_0, R$ (with $U:=T_{n,m,q}$). If (ii) of Lemma~\ref{thm-2-handle1} holds, then (ii) of this lemma holds. Now assume that (i) of Lemma~\ref{thm-2-handle1} holds. So the edge of $Q_1$ incident with $T_n$ should be a $(T_n,\phiv/R, \alpha_0, \dd_n)$-nonexit edge; but $Q_1$ is a $(T_n,\phiv/R,\alpha_0, \dd_n)$-exit path, a contradiction.  

Note that the most time-consuming steps above are the applications of Lemmas~\ref{thm-2-T(n,m)-exit} and \ref{lem-7-Pt-disjoint}. So, this proof takes $O(|V(G)|^4|E(G)|^2)$ time.   \qed

Next, we deal with the case when $z_\alpha\in V(branch(T)-v_0)$.

\begin{lemma} \label{lem-7-brach-cut} 
In $O(|V(G)|^5|E(G)|^2)$ time, 
{\em (i)} one concludes that $z_{\alpha}\notin V(branch(T)-v_0)$, or {\em (ii)} one finds a tree $T'$ and a coloring $\phiv'$ Kempe equivalent to $\phiv$, such that $T'$ is not $\phiv'$-elementary and $(s(T'),\ell(T'),p(T'),t(T'), b(T'))<(s(T),\ell(T),p(T),t(T), b(T))$.
\end{lemma}

\pf Suppose (i) fails and let $z_{\alpha}=v_i$ for some $i\in [t-1]$.  

First, suppose $i=t-1$. Let $R:=G[f_t]$. Note that $R\cap T_{n,m,q}=\emptyset$, and $\pi:=\phiv/R$ is Kempe equivalent to $\phiv$. Observe that if $R$ intersects a $(T_n, \phiv, \alpha_0, \dd_n)$-exit path then $R$ is edge disjoint from all $(T_n, \pi,\alpha_0, \dd_n)$-paths. Hence, 
the $(T_n, \phiv,\alpha_0, \dd_n)$-nonexit paths are precisely the $(T_n, \pi,\alpha_0, \dd_n)$-nonexit paths (or (ii) holds by Lemma~\ref{thm-2-Tn-extension}). 
Therefore, under $\pi$, $T^-$ is an $(n,m,q+1)$-Tashkinov tree.  Since $\phiv(f_t) \in \phibar(T^-)$, we have $\phiv(f_t)\in \pibar(v_{t-1})\cap \pibar(T(v_{t-1})^-)$. So (ii) holds with $T':=T^{-}$ and $\phiv':=\pi$.     

Now assume $i<t-1$.  We will increase $i$ to $i+1$ through Kempe exchanges, eventually reducing the proof to the case $i=t-1$. Let $\beta\in \phibar(v_{i+1})$. Since $T^-$ is $\phiv$-elementary, we have $\alpha,\beta\notin \phibar(T_{n, m, q})$; so $\alpha, \beta\notin \Gamma_q$.

Suppose $\{\alpha, \beta\}\cap S_n =\emptyset$. Then $\alpha, \beta\notin \{\alpha_0, \dd_n\}$. Let $R:= P_{v_t}(\alpha, \beta; \phiv)$ and $\pi=\phiv/R$. 
Apply Lemma~\ref{lem-7-Pt-disjoint} with $\epv_1=\alpha$, $\epv_2=\beta$, $v=v_i$, and $w=v_{i+1}$. We may assume $P_{v_i}(\alpha, \beta; \phiv) = P_{v_{i+1}}(\alpha, \beta; \phiv)$ (or (ii) holds by (i4) of Lemma~\ref{lem-7-Pt-disjoint}); so $v_i, v_{i+1}\notin V(R)$. We may also assume that $T$ is an $(n,m,q+1)$-Tashkinov tree under $\pi$ (or (ii) holds by (i5) of Lemma~\ref{lem-7-Pt-disjoint}). Now  $\beta\in \pibar(v_t)\cap \pibar(v_{i+1})$, increasing $i$ to $i+1$.

Therefore, let $\dd\in \{\alpha, \beta\}\cap S_n$. Then $\dd\in S_n\setminus \phibar(T_{n,m,q})$. Since $|\Gamma_q(\dd)| =2$,  there exists $\gamma\in \Gamma_q(\dd)\backslash\{\phiv(f_{i+1})\}$.  By definition, $\gamma\in (\phibar(T_{n,m,q})\setminus \{\alpha_0\})\backslash \phiv(T(v_{i+1}) - E(T_{n,m,q}))$. 
Let $R:= P_{v_t}(\alpha, \gamma; \phiv)$ and $\pi:= \phiv/R$. 
Apply Lemma~\ref{lem-7-Pt-disjoint} with $\epv_1 = \gamma$, $\epv_2=\alpha$, $v=z_{\gamma}$, and $w=v_{i}$. We may assume $P_{v_\gamma}(\alpha, \gamma; \phiv) = P_{v_i}(\alpha, \gamma; \phiv)$ (or (ii) holds by (i4) of Lemma~\ref{lem-7-Pt-disjoint}); so $v_\gamma, v_i\notin V(R)$. 
We may also assume that $T$ is an $(n,m,q+1)$-Tashkinov tree under $\pi$, or (ii) holds by (i5) of Lemma~\ref{lem-7-Pt-disjoint}. 

Note that $\gamma\in \pibar(T_{n,m,q})\cap \pibar(v_t)$. Let $R:=P_{v_t}(\gamma,\beta;\pi)$ and $\xi=\pi/R$. 
 Apply Lemma~\ref{lem-7-Pt-disjoint} again with $\epv_1 = \gamma$, $\epv_2 = \beta$, $v=z_{\gamma}$ and $w=v_{i+1}$. We may assume $P_{v_\gamma}(\beta, \gamma; \phiv) = P_{v_{i+1}}(\beta, \gamma; \phiv)$ (or (ii) holds by (i4) of Lemma~\ref{lem-7-Pt-disjoint}); so $v_\gamma, v_{i+1}\notin V(R)$.  
 We may also assume that $T$ is an $(n,m,q+1)$-Tashkinov tree under $\xi$, or (ii) holds by (i5) of Lemma~\ref{lem-7-Pt-disjoint}. Note that $\beta\in \xibar(v_t)\cap \xibar(v_{i+1})$, increasing $i$ to $i+1$. 

The above process takes $O(|V(G)|^4|E(G)|^2)$ time. By repeating it for $i<t-1$ at most $t<|V(G)|$ steps, we reach the case when $i=t-1$; hence (ii) holds from the argument in the second paragraph of the proof. Thus, the whole argument runs in $O(|V(G)|^5|E(G)|^2)$ time.\qed

\begin{lemma}\label{lem-7-s>0}
In $O(|V(G)|^5|E(G)|^2)$ time, 
{\em (i)} one concludes that $T-T_{n,m, q}$ is not a path
or  {\em (ii)} one finds a tree $T'$ and a coloring $\phiv'$ Kempe equivalent to $\phiv$, such that $T'$ is not $\phiv'$-elementary and $(s(T'),\ell(T'),p(T'),t(T'), b(T'))<(s(T),\ell(T),p(T), t(T), b(T))$.   
\end{lemma}

\proof We first comment that it will be clear from the proof that the complexity bound $O(|V(G)|^5|E(G)|^2)$ is due to the complexity in Lemma~\ref{lem-7-brach-cut}. Now suppose $T-T_{n,m,q}$ is a path (thus contained in $branch(T)$). Then we may assume $t\ge 2$ (or (ii) holds by Lemma~\ref{lem-7-1st-edge-elem}) and $z_{\alpha}\in V(T_{n,m, q})$ (or (ii) holds by Lemma~\ref{lem-7-brach-cut}). 

We may assume $\alpha\ne \alpha_0$. Otherwise, let $\alpha'\in \phibar(g)\backslash\{\alpha_0\}$ and $R:= P_{v_t}(\alpha',\alpha; \phiv)$ and $\pi=\phiv/R$. Then $R\cap T_{n,m,q} = \emptyset$ as $\alpha,\alpha'\notin \phiv(\partial T_{n,m,q})$.  
Apply Lemma~\ref{lem-7-Pt-disjoint} with $\{\epv_1,\epv_2\}=\{\alpha',\alpha\}$, $\{v,w\}=\{z_{\alpha'},z_{\alpha}\}$, and $v\in V(T(w))$. We may assume that $T_{n,m,q}$ is an $(n,m,q)$-tree under $\pi$ with the same sets of parameters as those under $\phiv$. Clearly, $T$ is a prefix of a closure of $T_{n,m,q}+g_{n,m,q}$. Since $\alpha=\alpha_0, \alpha'\in \phibar(G[g])$, it follows that $\alpha, \alpha'\notin \Gamma_q$, consequently, $T$ is an $(n,m,q)$-Tashkinov tree under $\pi$. 
  Hence, we could work with $\alpha',\pi$ instead of $\alpha,\phiv$.

{\it Case} 1. $\phiv(f_1)\ne \alpha$. 

Let $i=t-1$ if $\alpha \notin \phiv(branch(T))$; otherwise, let $i\in [t-1]$ be minimum such that $\phiv(f_{i+1}) = \alpha$. Then $\alpha\notin \phiv(T(v_i) - E(T_{n,m,q}))$ since $\alpha\in \phibar(T_{n,m,q})\backslash\{\alpha_0\}$. Let $\beta\in \phibar(v_i)$, $R:= P_{v_t}(\alpha, \beta; \phiv)$, and
$\pi = \phiv/R$. 

Apply Lemma~\ref{lem-7-Pt-disjoint} with $\epv_1=\alpha$, $\epv_2=\beta$, $v=z_{\alpha}$, and $w=v_i$. 
We may assume $P_{z_\alpha}(\alpha, \beta; \phiv) = P_{v_i}(\alpha, \beta; \phiv)$  (or (ii) holds by (i4) of Lemma~\ref{lem-7-Pt-disjoint}); so $z_\alpha, v_i\notin V(R)$. Note that if $i=t-1$, then $\phiv(f_t) \notin \{\alpha, \dd\}$; 
If $i < t$ and $\alpha\in \Gamma_q(\dd)$ for some $\dd\in S_n$, then $\dd\in \phibar(T(v_i))$ since $\phiv(f_{i+1}) =\alpha$.  Hence, under $\pi$, we may assume that $T$ is an  $(n,m,q+1)$-Tashkinov tree or (ii) holds by (i5) of Lemma~\ref{lem-7-Pt-disjoint}.  Since  $\beta\in \pibar(v_t)\cap \pibar(v_i)$, we see that (ii) follows from Lemma~\ref{lem-7-brach-cut}. 



{\it Case} 2. $\phiv(f_1)=\alpha$.  

If $q=0$, then
$\alpha=\dd_n \notin \phibar(T_{n,m,0})$ and hence, $z_{\alpha}\in V(branch(T)\backslash\{v_t\})$; then (ii) holds by Lemma~\ref{lem-7-brach-cut}.  
Thus, we may assume $q\ge 1$;  hence $f_1=g_{n,m,q}$ and  $\alpha\in \Gamma_{q-1}\setminus \Gamma_q$ (so $\alpha\notin \phibar(g)$).
Let $\alpha'\in \phibar(g)\backslash\{\alpha,\alpha_0\}$; then $\alpha'\notin \phiv(\partial T_{n,m,q})$. Let $R := P_{v_t}(\alpha, \alpha'; \phiv)$
and $\pi = \phiv/R$.  Apply Lemma~\ref{lem-7-Pt-disjoint} with $\epv_1 = \alpha'$, $\epv_2=\alpha$, $v=z_{\alpha'}$ and $w=z_{\alpha}$. We may assume that $P_{z_\alpha'}(\alpha, \alpha'; \phiv) = P_{z_\alpha}(\alpha, \alpha'; \phiv)$  (or (ii) holds by (i4) of Lemma~\ref{lem-7-Pt-disjoint}); so $z_\alpha, z_{\alpha'}\notin V(R)$. We may also assume that  $T_{n,m,q}$ is an $(n,m,q)$-Tashkinov tree under $\pi$, or (ii) holds by (i2) of Lemma~\ref{lem-7-Pt-disjoint}. 

Since $\alpha'\in \phibar(g)\backslash\{\alpha_0\}$, we may assume that there is at most one $(T_{n,m,q},\phiv,\alpha,\alpha')$-exit path; or or (ii) holds by Lemma~\ref{thm-2-T(n,m)-exit}. Thus, 
$R\cap T_{n,m,q} = \emptyset$, in particular, $f_1\notin E(R)$. 

Since  $\alpha, \alpha'\in \phibar(T_{n,m,q})\setminus \Gamma_q$, swapping $\alpha$ and $\alpha'$ on edges of $T-T_{n,m,q}$ does not affect the tree structure. Hence,  under $\pi$,  $T$ is an $(n,m,q)$-Tashkinov tree. Replacing $\phiv, z_\alpha, \alpha$ by $\pi,z_{\alpha'},\alpha'$, we get back to Case 1 where $\pi(f_1) =\phiv(f_1) \ne \alpha'\in \pibar(v_t)\cap \pibar (z_{\alpha'})$. \qed

\begin{lemma}\label{cla-4.7-1}
Suppose  $t \ge 2$ and $\eta \in \phibar(v_t)\cap \phibar(T(u_r)^-)\backslash\{\alpha_0\})$. Then 
 in $O(|V(G)|^5  |E(G)|^2)$ time, {\em (i)} one concludes that  $\eta\in \phiv(T(v_1) -E(T_{n,m,q}))$ or $\eta\in   S_n\cap \phibar(T(u_r)^- -T_{n,m,q})$, or {\em (ii)} one finds a tree $T'$ and a coloring $\phiv'$ Kemple equivalent to $\phiv$, such that $T'$ is not $\phiv'$-elementary and $(s(T'), \ell(T'), p(T'), t(T'), b(T')) < (s(T),\ell(T),p(T), t(T), b(T))$.
 \end{lemma}
 \pf First, we comment the complexity bound $O(|V(G)|^5  |E(G)|^2)$ is due to Lemma~\ref{lem-7-brach-cut}. Now turn to the proof, and suppose (i) fails; so  $\eta \notin \phiv(T(v_1) -E(T_{n,m,q}))$ (in particular $\eta\ne \phiv(f_1)$) and $\eta\notin S_n\cap \phibar(T(u_r)^- -T_{n,m,q})$.
 
 Let $i = t-1$ if $\eta\notin \phiv(branch(T) -v_0)$; otherwise, let $i\in [t]$ be minimum  such that $\phiv(f_{i+1})=\eta$.  Then $i \ge 1$ and  $\eta\notin \phiv(T(v_i) - E(T_{n,m,q}))$ since $\eta\notin \phiv(T(v_1)-E(T_{n,m,q}))$ and $\phiv(f_1) \ne \eta$, and we may assume that $T-T_{n,m,q}$ is not a path (or (ii) holds by Lemma~\ref{lem-7-s>0}). 
 Let $\dd\in \phibar(v_i)$. Then $\dd\notin \{\phiv(f_i), \eta\}$ and  $\dd\notin \phiv(T(v_i) - E(T_{n,m,q}))$
 unless $q=0$ and $\dd = \dd_n$, in this case, $\dd\notin \phiv(T(v_i) - E(T_n + g_{n,m}))$
 (as $\dd\notin \phibar(T(v_i)^-)$.   

{\it Case} 1.  $\eta\notin \phibar(T_{n,m,q})$ (so $z_\eta\notin V(T_{n,m,q})$) and $\dd\in S_n$  (so $\dd \in S_n\setminus \phibar(T_{n,m,q})$). 

 Let $\gamma\in \Gamma_q(\dd)$.  By definition, $\gamma\notin \phiv(T(v_i) - E(T_{n,m,q}))$. Apply Lemma~\ref{lem-7-Pt-disjoint} with $\epv_1=\gamma$, $\epv_2=\eta$, $v=z_{\gamma}$, $w=z_{\eta}$, $R := P_{v_t}(\eta, \gamma; \phiv)$,  and $\pi := \phiv/R$. We may assume $P_{z_\gamma}(\eta, \gamma; \phiv) = P_{z_\eta}(\eta, \gamma; \phiv)$ (or (ii) holds by (i4) of Lemma~\ref{lem-7-Pt-disjoint}); so $z_\gamma, z_\eta\notin V(R)$ and $R\cap T_{n, m, q} = \emptyset$ (or (ii) holds by Lemma~\ref{thm-2-Tn-extension}). Additionally, we notice that $\eta, \gamma\notin \phiv(T(v_i)-E(T(z_\eta)))$ and, hence, we may assume that $T(v_i)$ is an $(n,m,q+1)$-Tashkinov tree under $\pi$ (or (ii) holds by (i5) of Lemma~\ref{lem-7-Pt-disjoint}). 
 Note $\eta\notin \Gamma_q$ (as $\eta\notin \phibar(T_{n,m,q})$) and $\eta\in \phibar(T(u_r)^-)$ and $\gamma\in \Gamma_q(\dd)$ with  $\dd\in \phibar(v_i)$. So $T$ is obtained from $T(w)$ by following (b1) of Definition~\ref{defi-2-hi}; hence $T$ is an $(n,m,q+1)$-Tashkinov tree under $\pi$. 

 Note that $\gamma\in \pibar(v_t)\cap (\pibar(T_{n,m,q})\backslash\{\alpha_0\})$ and $\gamma\notin \pi(T(v_i) -E(T_{n,m,q}))$.  Apply Lemma~\ref{lem-7-Pt-disjoint} with $\epv_1 = \gamma$, $\epv_2 =\dd$, $v=z_\gamma$, $w=v_i$, $R :=P_{v_t}(\gamma, \dd; \pi)$, and $\xi = \pi/R$. We may assume $P_{z_\gamma}(\gamma, \dd; \pi) =P_{v_i}(\gamma, \dd; \pi)$ (or (ii) holds by (i4) of Lemma~\ref{lem-7-Pt-disjoint}); so $z_\gamma, v_i\notin V(R)$.  
 Since $v_i\notin V(T_{n,m,q}\cup R)$ and $\gamma\in \Gamma_q(\dd)$ and $\dd\in \phibar(v_i)$, we may assume that $T$ is an $(n,m,q+1)$-Tashkinov tree under $\xi$ (or (ii) holds by (i5) of Lemma~\ref{lem-7-Pt-disjoint}). 
 

{\it Case} 2. $\eta\in \phibar(T_{n,m,q})$ or $\dd\notin S_n$. 

By assumption, $\eta\ne\alpha_0$ and $\eta\notin S_n\cap \phibar(T(u_r)^- - T_{n,m,q})$. So, either $\eta\in \phibar(T_{n,m,q})\backslash\{\alpha_0\}$, or $\eta\notin \phibar(T_{n,m,q})$ in which case $\eta\notin S_n$
(as $\eta\in \phibar(T(u_r)^-)$ and $\dd\notin S_n$ and,  thus, $\{\eta,\dd\}\cap (\phibar(T_{n,m,q})\cup S_n)=\emptyset$.
Hence, we may apply Lemma~\ref{lem-7-Pt-disjoint} with $\epv_1=\eta$, $\epv_2 = \dd$, $v=z_\eta$, $w=v_i$,  $R:= P_{v_t}(\eta, \dd; \phiv)$ and $\pi=\phiv/R$. 
We may assume $P_{z_\eta}(\eta, \dd; \phiv) = P_{v_i}(\gamma, \dd; \phiv)$ (or (ii) holds by (i4) of Lemma~\ref{lem-7-Pt-disjoint});  so $z_\eta, v_i\notin V(R)$. 

If $\eta\notin \phiv(T- E(T(v_1)))$ then $i=t-1$ and so $\dd\notin \phiv(T-E(T(v_i)))$; otherwise, we have $\phiv(f_{i+1}) = \eta$, which in turn shows that if $\eta\in \Gamma_q(\dd)$ for some $\dd\in S_n$ then $\dd\in \phibar(T(v_i))$. Hence, the conditions of (i5) of Lemma~\ref{lem-7-Pt-disjoint} are satisfied. Thus, we may assume that $T$ is an $(n,m,q+1)$-Tashkinov tree under $\pi$, or else (ii) holds by (i5) of Lemma~\ref{lem-7-Pt-disjoint}. Now $\eta\in \pibar(v_t)\cap \pibar(v_i)$, and so (ii) holds by Lemma~\ref{lem-7-brach-cut}. \qed



\begin{lemma}\label{lem-7-swap-3}

 Let $\eta\in (\phibar(T(u_r)^-)\cap \phibar(v_t))\setminus \{\alpha_0\}$ with 
 $\eta\notin S_n\cap \phibar(T(u_r)^- -T_{n,m,q})$, and let 
 $\beta\in B(T(u_r)^-, \phiv) \backslash\{\phiv(f_1)\}$, 
 $R=P_{v_t}(\eta, \beta; \phiv)$, and $\pi= \phiv/R$. 
 Then in $O(|V(G)|^5|E(G)|^2)$ time,  
 \begin{itemize}
 \item [(i)] one concludes that {\em (ia)} $\eta\notin B(T(u_r)^-, \phiv)$ when $t\ge 2$, {\em (ib)} $\beta\in (B(T(u_r)^-, \pi)\cap \pibar(v_t))\setminus \{\pi(f_1)\}$,  {\em (ic)} 
 the $(T_n, \pi, \alpha_0, \dd_n)$-nonexit edges are precisely the $(T_n, \phiv, \alpha_0, \dd_n)$-nonexit edges, and {\em (id)} $T$ is an $(n,m,q+1)$-Tashkinov tree under $\pi$;  or 
 \item [(ii)] one finds a tree $T'$ and a coloring $\phiv'$ Kempe equivalent to $\phiv$, such that $T'$ is not $\phiv'$-elementary and  $(s(T'), \ell(T'), p(T'), t(T'), b(T')) < (s(T),\ell(T),p(T), t(T), b(T))$.
 \end{itemize}
 \end{lemma}

 \pf 
We comment first that the complexity bound is due to the application of 
Lemma~\ref{cla-4.7-1}. Now, turn to the proof. Since $\beta\in B(T(u_r)^-, \phiv)$, it follows that  $\beta\in \phibar(T(u_r)^-)\backslash\{\alpha_0\}$ and 
$\beta\notin S_n\cap (\phibar(T(u_r)^- \setminus E(T_{n,m,q})))$. 
Hence, by assumption on $\eta$, $\{\eta, \beta\}\cap (\phibar(T_{n,m,q})\backslash \{\alpha_0\}) \ne \emptyset$ or $\{\eta, \beta\}\cap (\phibar(T_{n,m,q}) \cup S_n) = \emptyset$. 
So we could apply Lemma~\ref{lem-7-Pt-disjoint} with $\{\epv_1,\epv_2\}=\{\eta,\beta\}$ in Case 1 below.

{\it Case} 1. $R\cap T(u_r)\ne \emptyset$.

 Let $Q_w$ be the subpath of $R$ from $v_t$ to $w\in V(T(u_r))$ and internally disjoint from $T(u_r)$. Let $T':=T(u_r)\cup P_w$ if $w=u_r$; and $T':=T(u_r)^-\cup Q_w$ if $w\ne u_r$. The edges of $T'$ are ordered as in $T(u_r)$ followed by the edges in $Q_w$ from $w$ to $v_t$. 
Then $t(T')<t(T)$, because if $w=u_r$ then $branch(T')$ contains $Q_w$ and the last edge of $T(u_r)$, and if $w\ne u_r$ then $trunk(T')$ is contained in $T(u_r)^-$. 

We may assume that $T'$ is not an $(n,m,q+1)$-Tashkinov tree under $\phiv$; as otherwise (ii) holds with $T'$ and $\phiv':=\phiv$. Therefore, $T'$ violates (b1) of Definition~\ref{defi-2-hi}:   $\eta$ or  $\beta \in \Gamma_q(\dd)$ for some $\dd\in S_n \setminus \phibar(T(u_r)^-)$. In particular, $\eta$ or $\beta \notin B(T(u_r)^-, \phiv)$, and so $\eta \notin B(T(u_r)^-, \phiv)$ (so (ia) holds), and either $\dd\notin \phibar(T^-)$ or $z_\dd\notin V(T(u_r)^-)$. In the former case, let $T(z_\dd)= T$. By definition,   $\eta\notin \phiv(T(z_\dd) -  E(T_{n,m,q}))$. 

{\it Subcase} 1.1 $z_\dd\notin \phibar(T(u_r))$.

Suppose $t=1$. Apply Lemma~\ref{lem-7-Pt-disjoint} with $\{\epv_1,\epv_2\}=\{\eta,\beta\}$, $\{v,w\}=\{z_\eta,z_\beta\}$, and $v\in V(T(w))$.  
We may assume that $P_{z_\eta}(\eta, \beta; \phiv) = P_{z_\beta}(\eta, \beta; \phiv)$ (or (ii) holds by (i4) of Lemma~\ref{lem-7-Pt-disjoint}); so $z_\eta, z_\beta\notin V(R)$. We may also assume that (ic) holds (or (ii) holds by (i1) of Lemma~\ref{lem-7-Pt-disjoint}) and $T_{n,m,q}$ is an $(n,m,q)$-tree (or (ii) holds by (i2) of Lemma~\ref{lem-7-Pt-disjoint}).   Since $\eta, \beta\notin \phiv(T(v_1) -  E(T_{n,m,q}))$, $T=T(v_1)$ is an $(n,m,q+1)$-Tashkinov tree under $\pi$ (so (id) holds). Clearly, $\beta\in B(T(u_r)^-,\pi)\cap \pibar(v_t)\setminus \{\pi(f_1)\}$ (so (i2) holds). Hence, we have (i). 

Now assume $t\ge 2$. 
Recall that $\eta \notin S_n\cap \phibar(T(u_r)^- -T_{n,m,q})$ and $\eta\notin \phiv(T(v_1) -E(T_{n,m,q}))$ since $z_\dd\notin V(T(u_r))$, which in turn shows that $\eta\notin \phibar(v_t)$, which gives a contradiction to   $\eta\in   \phibar(v_t)$, and so (ii) holds by Lemma~\ref{cla-4.7-1}.

{\it Subcase} 1.2. $z_\dd=u_r$; so $\dd\in \phibar(u_r)$. 

We now apply Lemma~\ref{lem-7-Pt-disjoint} with $\epv_1$ and $\epv_2$ as in the $t=1$ case of Subcase 1.1. We may assume that $z_\eta, z_\beta\notin V(R)$, else (ii) holds by (i4) of Lemma~\ref{lem-7-Pt-disjoint}. We may also assume that (ic) holds (else else (ii) holds by (i1) of Lemma~\ref{lem-7-Pt-disjoint}) and $T_{n,m,q}$ is an $(n,m, q)$-tree under $\pi$ (or (ii) holds by (i2) of Lemma~\ref{lem-7-Pt-disjoint}). 

Since $\eta, \beta\in \phibar(T(u_r)^-)$ and $z_\eta, z_\beta\notin R$, we have 
$\pi(f_i) \in \pibar(T(f_i))$ for $i\in [t]$.
Since $\eta, \beta\notin \phiv(T(u_r) - E(T_{n,m,q}))$, we see that  $T(u_r)$, under $\pi$,  is an $(n,m,q)$-Tashkinov tree. 
Note that $\beta\in \Gamma_q(\dd')$ implies $\dd'\in \phibar(T(u_r)^-)$ (as $\beta\in B(T(u_r)^-, \phiv)$) and  that $\eta\in \Gamma_q(\dd)$ with $\dd\in \phibar(u_r)$. Hence, by (b1) of Definition~\ref{defi-2-hi},    
$T$ is an $(n,m,q+1)$-Tashkinov tree under $\pi$ (so (id) holds). Since $\eta, \beta\notin \phiv(T(u_r))$, $\beta\in B(T(u_r)^-, \pi)$ (so (ib) holds). Hence, we have (i).

{\it Case} 2. $R\cap T(u_r)=\emptyset$. 

Then $\pi$ and $\phiv$ agree on $T(u_r)$, $\beta\in B(T(u_r)^-, \pi)\cap \pibar(v_t)\setminus \{\pi(f_1)\}$ (so (ib) holds), and $T_n$ is an $(n,0,0)$-tree under $\pi$. 

We claim that (ic) holds and $T(u_r)$, under $\pi$, is an $(n,m,q)$-Tashkinov tree. This is clear if $\{\eta, \beta\}\cap \{\alpha_0, \dd_n\} = \emptyset$. Now suppose $\{\eta, \beta\}\cap \{\alpha_0, \dd_n\} \ne \emptyset$. Then $\dd_n\in \{\eta, \beta\}$ as $\alpha_0\notin \{\eta, \beta\}$. Since $\eta, \beta\notin S_n\cap (\phibar(T(u_r)^- -T_{n,m,q}))$, we have $\dd_n\in \phibar(T_{n,m,q})$. By definition, $\dd_n\notin \Gamma_q$; so $\dd_n\notin \partial T_{n,m,q}$. Hence $V(T_n^+(\phiv, \eta, \beta))\subseteq V(T_{n,m,q})$ which is disjoint from $R$. Thus, (ic) holds, and $T(u_r)$ is an $(n,m,q+1)$-Tashkinov tree under $\pi$.

 Since both $\eta, \beta\in \phibar(T(u_r)^-)$, it follows that $T$ is a prefix of a closure of $T(u_r)$.  
 Since $\beta\in B(T(U_r)^-, \phiv)$, if $\beta \in \Gamma_q(\dd)$ then $\dd\in S_n\cap (\phibar(T(u_r)^--T_{n,m,q}))$. 
 By the assumption and Lemma~\ref{cla-4.7-1}, we have $\eta\in \phiv(T(v_1))$, consequently, 
 if  $\eta\in \Gamma_q(\dd)$ then $\dd\in S_n\cap (\phibar(T(u_r)^--T_{n,m,q}))$. Hence, (b1) of Definition~\ref{defi-2-hi} is satisfied, and so  $T$ is an $(n, m,q)$-Tashkinov tree. So (id) holds. 
 
If (ia) holds then we have (i). So, assume $t\ge 2$ and $\eta\in B(T(u_r)^-, \phiv)$. Since $R\cap T(u_r)=\emptyset$, we see that $\eta, \beta\in B(T(u_r)^-, \pi)$. Hence, $\eta,\beta\notin S_n\cap \phibar(T(u_r)^-)-T_{n,m,q})$. Since $\eta\in \phibar(v_t)$, we may assume that $\eta\in \phiv(T(v_1)- E(T_{n,m,q}))$ (or (ii) holds by Lemma~\ref{cla-4.7-1}). But then   $\beta\in \pibar(v_t)$. Since $\beta\in B(T(u_r)^-, \pi)$, we have $\beta\notin \pi(T(v_1) - E(T_{n,m,q}))$. Hence, (ii) holds by Lemma~\ref{cla-4.7-1}). \qed

\begin{lemma}\label{lem-7-trunk-clean}
In $O(|V(G)|^5|E(G)|^2)$ time,
\vspace{-2ex}
\begin{itemize}
\item [(i)] one concludes that $t=1$ and finds a coloring $\pi$ Kempe equivalent to $\phiv$ such that {\em (ia)} $T$, under $\pi$, is an $(n,m,q+1)$-Tashkinov tree, {\em (ib)} $\pibar(v_1)\cap B(T(u_r)^-, \pi) \not\subseteq \{\pi(f_1)\}$, and  if $\phiv(f_1)\notin \phibar(u_r)$ and $\phiv(f_1)\notin  \bigcup_{\dd\in S_n\cap \phibar(u_r)}\Gamma_q(\dd)$, then $\pi(f_1)\notin \pibar(u_r)$ and $\pi(f_1)\notin \bigcup_{\dd\in S_n\cap \pibar(u_r)}
 \Gamma_q(\dd)$; or 
\item [(ii)] one finds a tree $T'$ and a coloring $\phiv'$ Kempe equivalent to $\phiv$, such that $T'$ is not $\phiv'$-elementary and  $( s(T'), \ell(T'),p(T'),t(T'),b(T') ) <  (s(T), \ell(T),p(T), t(T), b(T) )$.

\end{itemize}
\end{lemma}

\pf First, we remark that the complexity bound $O(|V(G)|^5|E(G)|^2)$ is due to applications of the previous lemmas. Now for the proof, we may assume $u_r\notin V(T_{n,m,q})$ (or (ii) holds by Lemma~\ref{lem-7-s>0}) and $\alpha \in \phibar(v_t)\cap \phibar(T(u_r))$ (or (ii) holds by Lemma~\ref{lem-7-brach-cut}.

 Suppose $\alpha=\alpha_0$. 
Let $\alpha' \in \phibar(g)\backslash\{\alpha_0\}$. Apply Lemma~\ref{lem-7-Pt-disjoint} with $\{\epv_1,\epv_2\}=\{\alpha',\alpha\}$, $\{v,w\}=\{z_{\alpha'},z_\alpha\}$, $v\in V(T(w))$, and $R:=P_{v_t}(\alpha, \alpha'; \phiv)$. We may assume 
$z_{\alpha}, z_{\alpha'}\notin V(R)$ (or (ii) holds by (i4) of Lemma~\ref{lem-7-Pt-disjoint}) and $T_{n,m,q}$ is an $(n,m,q)$-tree. We may assume that  $R\cap T_{n,m,q} = \emptyset$ (or (ii) holds by Lemma~\ref{thm-2-T(n,m)-exit}).  Since $\{\alpha', \alpha_0\}\subset \phibar(G[g])\setminus \Gamma_q$, we may assume that $T$ is an $(n,m,q)$-Tashkinov tree by Definition~\ref{defi-2-hi}, and could consider $\phiv/R$ instead of $\phiv$.
Hence, we may assume $\alpha \ne \alpha_0$. We consider three cases.

{\it Case} 1.  $\alpha \in \phibar(v_t)\cap \phibar(u_r)$ and $\alpha\notin S_n\cap \phibar(T(u_r)-T_{n,m,q})$.

Then $\alpha\notin \phibar(T_{n,,m,q})$ (as $T^-$ is $\phiv$-elementary) and, hence, in this case we have $\alpha\notin \phibar(T_{n,m,q})\cup S_n$ and $\alpha\notin \phiv(T(u_r) - E(T_{n,m,q}))$.
Since $|B(T_{u_r}^-, \phiv)| \ge 2$, there exists $\eta\in B(T_{u_r}^-, \phiv)\backslash\{\phiv(f_1)\}$. In particular, $\eta\in \phibar(T_{n,m,q})$ or $\eta\notin \phibar(T_{n,m,q})\cup S_n$, and $\eta\notin \phiv(T(u_r) - E(T_{n,m,q}))$. Therefore,  $\{\alpha, \eta\}\cap (S_n\cup \phibar(T_{n,m,q}))=\emptyset$ if $\eta\notin \phibar(T_{n,m,q})\backslash\{\alpha_0\}$. Now apply Lemma~\ref{lem-7-Pt-disjoint} with $\epv_1=\eta$, $\epv_2=\alpha$, $v=z_\eta$, $w=u_r$, $R=P_{v_t}(\alpha, \eta; \phiv)$, and $\pi:= \phiv/R$.
By (i4), we may assume that $z_\eta, u_r\notin V(R)$ (else (ii) holds by (i4) of Lemma~\ref{lem-7-Pt-disjoint}) and that $T$ is an $(n,m,q+1)$-Tashkinov tree under $\pi$ (else (ii) holds by (i5) of Lemma~\ref{lem-7-Pt-disjoint}, since $u_r\notin V(T_{n,m,q}\cup R)$ and if $\eta\notin \Gamma_q(\dd)$ then $\dd\in \phibar(T(u_r)^-)$).  Now  $\eta\in \pibar(v_t)\cap \pibar(T(u_r)^-)$.

Since $\alpha\in \phibar(u_r)\backslash S_n$, $\alpha\notin \phiv(T(u_r) - E(T_n))$ and, hence, $\eta \notin \pi(T(u_r) \setminus E(T_n))$.  Thus, since
$\eta\in B(T(u_r)^-, \phiv)$, we have $\eta\in B(T(u_r)^-, \pi)$. Therefore, we can apply Lemma~\ref{lem-7-swap-3} to $T, \pi,\eta$ and some $\beta\in B(T(u_r)^-, \pi)\setminus \{\pi(f_1)\}$; so we may  assume $t=1$, else (ii) holds by (ia) of Lemma~\ref{lem-7-swap-3}.  Then $\eta\notin \pi(f_1)$ as $\eta\in \pibar(v_1)$.

Now suppose $\phiv(f_1)\notin \phibar(u_r)$ and  $\phiv(f_1) \notin \Gamma_q(\dd))$ for any $\dd \in S_n\cap \phibar(u_r)$. Then $\phiv(f_1) \notin \{\alpha,\eta\}$; so $f_1\notin E(R)$ and $\pi(f_1) = \phiv(f_1)$.   
If $u_r\in V(R)$ then we see that $T':=T(u_r)\cup R$ is an $(n,m,q)$-Tashkinov tree with $t(T')<t(T)$ (since $\alpha\in \phibar(u_r)$ and if $\eta\in \Gamma_q(\dd^*)$ then $\dd^*\in \phibar(T(u_r)^-)$); so (ii) holds. Hence, we may assume that $u_r\notin V(R)$. Then  $\pibar(u_r) = \phibar(u_r)$. So $\pi(f_1)\notin \pibar(u_r)$ and $\pi(f_1)\notin \Gamma_q(\dd)$ for $\dd\in S_n\cap \pibar(u_r)$.

{\it Case} 2. $\alpha\in \phibar(v_t)\cap \phibar(T_{u_r}^-)$ and $\alpha\notin S_n\cap \phibar(T(u_r)-T_{n,m,q})$. 

Choose $\beta\in B(T_{u_r}^-, \phiv)\backslash\{\phiv(f_1)\}$. We apply Lemma~\ref{lem-7-swap-3} with $\eta:=\alpha$, 
$R:=P_{v_t}(\alpha, \beta; \phiv)$, and $\pi=\phiv/R$.
 We may assume that $T$ is an $(n,m,q+1)$-Tashkinov tree under $\pi$ and  $\beta\in \phibar(v_t)\cap B(T(u_r)^-, \pi)$ (or (ii) holds by (id) and (ic) of Lemma ~\ref{lem-7-swap-3}). 
 
 We further assume that $t=1$;  as otherwise (ii) holds by  Lemma~\ref{cla-4.7-1}  (with $\eta:=\beta$). Now $\pi(f_1)\notin \{\alpha,\eta\}$ as $\{\alpha,\beta\}\subseteq \phibar(v_1)$.

Suppose $\phiv(f_1)\notin \phibar(u_r)$ and $\phiv(f_1) \notin \Gamma_q(\dd)$ for any $\dd\in S_n\cap \phibar(u_r)$. 
Since $\alpha, \eta\notin \phibar(u_r)$, $R$ does not end at $u_r$; so $\pibar(u_r) = \phibar(u_r)$.  
Since $f_1\notin E(R)$,  
$\pi(f_1)=\phiv(f_1)$.  Hence,  $\pi(f_1) \notin \pibar(u_r)$ and $\pi(f_1)\notin \Gamma_q(\dd)$ for $\dd\in S_n\cap \phibar(u_r)$.

{\it Case} 3. $\alpha \in S_n\cap \phibar(T(u_r) - T_{n,m,q})$.

Then $z_\alpha\notin V(T_{n,m,q})$. We will reduce this case to Case 2. Since $|\Gamma_q(\alpha)| =2$, there exists $\gamma\in \Gamma_q(\alpha)\backslash\{\phiv(f_1)\}$. Then $\gamma \in \phibar(T_{n,m,q})\backslash\{\alpha_0\}$ and $\gamma \notin \phiv(T(z_\alpha)-E(T_{n,m,q}))$.
We apply Lemma~\ref{lem-7-Pt-disjoint} with $\epv_1=\gamma$, $\epv_2=\alpha$, $v=z_\gamma$, $w=z_\alpha$,  $R:=P_{v_t}(\alpha, \gamma; \phiv)$, and $\pi=\phiv/R$. We may assume that $z_\gamma, z_\alpha\notin V(R)$ (else (ii) holds by (i4) of Lemma~\ref{lem-7-Pt-disjoint})  and that $T(z_\alpha)$ is an $(n,m,q+1)$-Tashkinov tree under $\pi$ (else (ii) holds by (i3) of Lemma~\ref{lem-7-Pt-disjoint}). Hence, $T$ is an $(n,m,q+1)$-Tashkinov tree,  
since $\alpha\notin \Gamma_q$ and $\gamma\in \Gamma_q(\alpha)$.


Suppose $\phiv(f_1)\notin \phibar(u_r)$ and $\phiv(f_1)\notin \Gamma_q(\dd)$ for $\dd\in S_n\cap \phibar(u_r)$. 
Since $\alpha\ne \phiv(f_1)$ and $\gamma\ne \phiv(f_1)$ (by choice), we have $f_1\notin E(R)$. 
Since $\{\alpha, \gamma\}\cap \phibar(u_r) = \emptyset$,  $R$ does not end at $u_r$; so $\pibar(u_r) = \phibar(u_r)$. Note that $\Gamma_q(\dd)$ remains the same under $\phiv$ and $\pi$. Hence, $\pi(f_1)=\phiv(f_1)\notin \pibar(u_r)$ and $\pi(f_1)\notin \Gamma_q(\dd)$ for any $\dd\in S_n\cap \pibar(u_r)$. 

Moreover, note that $\gamma\notin S_n\cap \pibar(T(u_r) - T_{n,m,q})$. 
Hence, by considering $\gamma, \pi$ instead of $\alpha, \phiv$, we get back to Case 2. \qed

We can now prove the main lemma needed to prove Theorem~\ref{thm:main}.

\begin{lemma}\label{lem-final}
 In $O(|V(G)|^5|E(G)|^2)$ time, one finds a Tashkinov tree $T'$ and coloring $\pi$ Kempe equivalent to $\phiv$, such that $T'$ is not $\phiv'$-elementary and\\ $(s(T'),\ell(T'),p(T'),t(T'),b(T'))< (s(T), \ell(T),p(T),t(T),b(T))$.
\end{lemma}
\pf  We may assume $t=1$ and 
$\alpha\in B(T(u_r)^-,\phiv)\setminus \{\phiv(f_1)\}$; otherwise, the assertion follows from Lemma~\ref{lem-7-trunk-clean}. Let $\theta=\phiv(f_1)$. Note $\theta\ne \alpha$, as $\alpha\in \phibar(v_t)$.  We may also assume that $u_r\notin V(T_{n,m,q})$,
otherwise the assertion holds by Lemma~\ref{lem-7-s>0}.

We consider three cases and reduce Case 2 and Case 3 to Case 1. We remark that the complexity $O(|V(G)|^5|E(G)|^2)$ is due to the application of lemmas in this section.

{\it Case} 1. $\phiv$ can be chosen so that  $\theta\notin \phibar(u_r)$ and $\theta\notin \Gamma_q(\dd)$ for any $\dd\in S_n\cap \phibar(u_r)$.

Let $T'$ denote the tree $T-u_r$ inheriting the edge ordering of $T$.  Then $T'$ is an $(n,m,q)$-Tashkinov tree under $\phiv$. Note that $(t(T'),b(T'))<(t(T),b(T))$; so the assertion holds with $T'$ and $\phiv':=\phiv$.

{\it Case} 2.  $\theta\in \phibar(u_r)$ or $\theta\in \Gamma_q(\dd)$ for some $\dd\in S_n\cap \phibar(u_r)$; and  $\theta\notin S_n\backslash\phibar(T_{n,m,q})$.

Suppose $\theta\in\phibar(u_r)$. Then $\theta\notin \phibar(T(u_r)^-)$ as $T^-$ is $\phiv$-elementary; so $\theta\notin \phibar(T_{n,m,q})\cup S_n$ (as $\theta\notin S_n\backslash\phibar(T_{n,m,q})$).   Since $\alpha\in B(T(u_r)^-, \phiv)$, we have $\alpha \ne \alpha_0$, $\alpha\notin \phiv(T(u_r)- E(T_n))$, and $\alpha\notin S_n\setminus  \phibar(T_{n,m,q})$. Therefore, we may apply Lemma~\ref{lem-7-Pt-disjoint} with $\epv_1 = \alpha$,  $\epv_2=\theta$, $v=z_\alpha$, $w=u_r$, $R:= P_{v_1}(\alpha,\theta;\phiv)$, and $\phiv_1 = \phiv/R$. We may assume $z_\alpha, u_r\notin V(R)$ (else the assertion follows from (i4) of Lemma~\ref{lem-7-Pt-disjoint}); so $\phibar_1(u_r)=\phibar(u_r)$.  We may also assume that $T$ is an $(n,m,q+1)$-Tashkinov tree under $\phiv_1$, else the assertion follows from (i5) of Lemma~\ref{lem-7-Pt-disjoint} (as $\dd'\in S_n\cap \phibar(T(u_r)^-)$ if $\alpha\in \Gamma_q(\dd)$). 
Now $\phiv_1(f_1)=\alpha\notin \phibar_1(u_r)$ and $\phiv_1(f_1)\notin \Gamma_q(\dd')$ for $\dd'\in S_n\cap \phibar_1(u_r)$. We apply Lemma~\ref{lem-7-trunk-clean} to $\phiv_1$. Note that (ii) of Lemma~\ref{lem-7-trunk-clean} gives the assertion of this lemma; so assume  (i) of Lemma~\ref{lem-7-trunk-clean} holds. Then there is a coloring $\xi$ Kempe equivalent to $\phiv_1$ such that $T$, under $\xi$, is an $(n,m,q+1)$-Tashkinov tree,
$\xibar(v_t)\cap B(T(u_r)^-, \xi) \not\subseteq \{\xi(f_1)\}$, 
 $\xi(f_1)\notin \xibar(u_r)$, and $\xi(f_1)\notin \Gamma_q(\dd')$ for $\dd'\in S_n\cap \xibar(u_r)$;  and we get back to Case 1.

Now assume that $\theta\in \Gamma_q(\dd)$ for some $\dd\in S_n\cap \phibar(u_r)$. 
Let $\gamma\in \Gamma_q(\dd)\setminus \{\theta\}$. Note that $\gamma\in \phibar(T_{n,m,q})\setminus \{\alpha_0\}$ and $\gamma\notin \phiv(T(u_r)-E(T_{n,m,q}))$. So we apply Lemma~\ref{lem-7-Pt-disjoint} with $\{\epv_1,\epv_2\}=\{\alpha,\gamma\}$, $\{v,w\}=\{z_\alpha,z_\gamma\}$ with $v\in V(T(w))$, $R:= P_{v_1}(\alpha, \gamma; \phiv)$, and $\phiv_1 = \phiv/R$.  Then $z_\alpha, z_\gamma\notin V(R)$ (else the assertion holds by (i4) of Lemma~\ref{lem-7-Pt-disjoint}); so $\phibar_1(u_r)=\phibar(u_r)$. Also,   
$T_{n,m,q}$ is an $(n,m,q)$-tree under $\phiv_1$ (else the assertion holds by (i3) of Lemma~\ref{lem-7-Pt-disjoint});  $T(u_r)$ is an $(n,m,q+1)$-Tashkinov tree under $\phiv_1$ as $\alpha, \gamma\in \phibar(T(u_r)^-)$ and $\alpha, \gamma\notin \phiv(T(u_r)-E(T_{n,m,q}))$. 
Note that $\phiv_1(f_1)=\alpha\notin \phibar_1(u_r)$ and  $\alpha\notin \Gamma_q(\dd')$ for $\dd'\in S_n\cap \phibar_1(u_r)$ (because $\alpha\in B(T(u_r)^-, \phiv)$).  Hence, $T$ is an $(n,m,q+1)$-Tashkinov tree under $\phiv_1$, and we can apply Lemma~\ref{lem-7-trunk-clean} in the same way as the previous paragraph and get back to Case 1 with $\xi$ replacing $\phiv$.

{\it Case} 3. $\theta\in \phibar(u_r)$ or $\theta\in \Gamma_q(\dd)$ for some $\dd\in S_n\cap \phibar(u_r)$; but $\theta\in S_n\backslash\phibar(T_{n,m,q})$.

Since $\theta\in S_n\backslash\phibar(T_{n,m,q})$, we have $\theta\notin \Gamma_q$, and so $\theta\in \phibar(u_r)$.  Let $\gamma_1, \gamma_2\in \Gamma_q(\theta)$. 

Let  $\phiv_1 = \phiv/P_{v_1}(\alpha,\gamma_1; \phiv)$. Note that $\gamma_1\in \phibar(T(u_r)^-)\backslash\{\alpha_0\}$, $\gamma_1\notin S_n\cap \phibar(T(u_r) - T_{n,m,q})$, and $\alpha\in B(T(u_r)^-, \phiv)$. 
 We may assume that $\gamma_1\in B(T(u_r)^-,\phiv_1)$,  
 $(T_n, \phiv_1, \alpha_0, \dd_n)$-nonexit edges are precisely the $(T_n, \phiv, \alpha_0, \dd_n)$-nonexit edges, and $T$, under $\phiv_1$, is an $(n,m, q+1)$-Tashkinov tree; otherwise the assertion holds by Lemma~\ref{lem-7-swap-3}.
 Note that $\phiv_1(f_1) =\theta$ and $\gamma_1\in \phibar_1(v_1)$.

  Since $\gamma_1\in \phibar(T_{n,m,q})$ and $\gamma_1\notin \phiv(T(u_r)-E(T_{n,m,q}))$, we apply Lemma~\ref{lem-7-Pt-disjoint} with $\epv_1= \gamma_1$, $\epv_2=\theta$, $v=z_{\gamma_1}$, $w=u_r$, $R:=P_{v_1}(\gamma_1, \theta; \phiv_1)$, and $\phiv_2 = \phiv_1/R$.  Then $z_\gamma, u_r\notin V(R)$ (else the assertion holds by (i4) of Lemma~\ref{lem-7-Pt-disjoint}) and $T$ is an $(n,m,q+1)$-Tashkinov tree under $\phiv_2$ (else the assertion holds by (i5) of Lemma~\ref{lem-7-Pt-disjoint}). Note that $\phiv_2(f_1) =\gamma_1\notin \phibar(u_r)$ and $\phiv_2(f_1)\notin \Gamma_q(\dd')$ for $\dd'\in S_n\cap \phibar(u_r)$. 

Apply Lemma~\ref{lem-7-trunk-clean} to $T, \phiv_2$ and $\theta\in \phibar_2(v_1)\cap \phibar_2(u_r)$. We may assume that there is a coloring $\xi$ Kempe equivalent to $\phiv_2$ (hence $\phiv$) such that $T$, under $\xi$, is an $(n,m,q+1)$-Tashkinov tree,
$\xibar(v_t)\cap B(T(u_r)^-, \xi) \not\subseteq \{\xi(f_1)\}$, 
 $\xi(f_1)\notin \xibar(u_r)$, and $\xi(f_1)\notin \Gamma_q(\dd')$ for $\dd'\in S_n\cap \xibar(u_r)$. Again, we get back to Case 1 with $\xi$ replacing $\phiv$. \qed  


\section{Conclusion}

We now complete the proof of Theorem~\ref{thm:main}. Let $G$ be a multigraph, $k:=\max\{\Delta(G)+1, \Gamma(G)\}$.

{\bf Step 0}.  Greedily obtain partial coloring $\phiv$ of $G$ with $\phiv(G)\subseteq [k]$. If $\phiv^{-1}([k])=E(G)$, we are done and output $\phiv$. Otherwise, choose $g\in E(G)\setminus \phiv^{-1}([k])$ and go to Step 1. This step takes $O(|V(G)|k)$ time. 

{\bf Step 1}. We grow an $(n,0,0)$-tree $T_{n}$, with base $T_0=cl_{\phiv}(g)$, stages $T_1, \ldots, T_n$, connecting colors $S_n:=\{\dd_0,\delta_1, \ldots, \delta_n\}$, such that $T_{i}=cl_{\phiv}(T_{i-1}^+(\phiv, \alpha_0,\dd_{i-1}))$ for $i\in [n]$, and $T_n$ is $\phiv$-elementary. We maximize $n$. This step takes $O(|V(G)|^2k^2)$ time. 

{\bf Step 2}. By the maximality of $n$, $T_{n+1}=cl_{\phiv}((T_n^+,\phiv, \alpha_0,\dd_{n+1}))$ is not $\phiv$-elementary. We grow $T_{n,0}:=T_n$ to an $(n,m,0)$-tree by adding levels $T_{n,1}, \ldots, T_{n,m}$ with extending edges $g_{n,0}, \ldots, g_{n,m-1}$ such that $T_{n,i}=cl_{\phiv}(T_{n,i-1}+g_{n,i-1},\phiv,\alpha_0,\dd_{n})$ for $i\in [m]$,  and $T_{n,m}$ is $\phiv$-elementary.  We  maximize $m$. This step takes $O(|V(G)|^2k^2)$ time. 

{\bf Step 3}.
By the maximality of $m$, $V(T_{n,m+1})\subseteq V(T_{n+1})$ and $T_{n,m+1}$ is not $\phiv$-elementary. We grow $T_{n,m,0}:=T_{n,m}$ by adding phases $T_{n,m,1}, \ldots, T_{n,m,q}$ such that $T_{n,m,q}$ is $\phiv$-elementary. We maximize $q$. 
This step takes $O(|V(G)|^2k^2)$ time.

{\bf Step 4}. By the maximality of $q$,  $V(T_{n,m,q+1})\subseteq V(T_{n,m+1})$ and $T_{n,m,q+1}$ is not $\phiv$-elementary. Choose the prefix $T$ of $T_{n,m,q+1}$ such that $T$ is not $\phiv$-elementary and, subject to this, $T$ is minimal. Then $T$ is an $(n,m,q)$-Tashkinov tree. If $n\ne 0$, go to Step 5. Now suppose $n=0$ (i.e. $T$ has no stages). Apply Lemma~\ref{lem:recoloring} and obtain a coloring $\pi$ Kempe equivalent to $\phiv$ such that $\pi^{-1}([k])=\phiv^{-1}([k])\cup \{g\}$. If $\pi^{-1}([k])=E(G)$, stop and output $\pi$; otherwise, update $\phiv$ with $\pi$, update $g$ with an edge from $E(G)\setminus \pi^{-1}([k])$,  go to Step 1. This step takes $O(|V(G)|^3|E(G)|)$ time.

{\bf Step 5}. Apply Lemma~\ref{lem-final}.  We obtain an ordered tree $T'$ that is not $\pi$-elementary under a coloring $\pi$ Kempe equivalent to $\phiv$, such that 
$(s(T'),\ell(T'),p(T'),t(T'), b(T'))< (s(T), \ell(T), p(T), t(T), b(T))$. Update $T$ with $T'$ and $\phiv$ with $\pi$, and go to Step 4. This step takes $O(|V(G)|^5|E(G)|^2)$ time.
 
Since the entire coloring process repeats Steps 1--5  at most $|E(G)|$ times and $k=O(|E(G)|)$, we arrive at the desired coloring in 
$O(|V(G)|^5|E(G)|^3)$ time.

\bibliographystyle{plain}
\bibliography{fcoloring.bib}

\end{document}